\documentclass[11pt]{amsart}

\usepackage[all]{xy}
\usepackage{latexsym,calc,amssymb,dsfont,amsfonts}
\usepackage{amsmath,amsthm,epsfig,tikz-cd}
\usepackage{graphicx,caption,subcaption}%,hyperref}
\usepackage[includefoot=true]{geometry}
\usepackage[utf8]{inputenc}
\usepackage{hyperref}

\newtheorem{thm}{Theorem}
\newtheorem{cor}[thm]{Corollary}
\newtheorem{lem}[thm]{Lemma}
\newtheorem{prop}[thm]{Proposition}
\theoremstyle{definition}
\newtheorem{defn}[thm]{Definition}
\theoremstyle{remark}

\newtheorem{problem}{Problem}

\newtheorem{examplex}{Example}
\newenvironment{example}
  {\pushQED{\qed}\examplex}
  {\popQED\endexamplex}

\DeclareMathAlphabet{\mymathbb}{U}{BOONDOX-ds}{m}{n}

\newcommand{\overbar}[1]{\mkern 2mu\overline{\mkern-2mu#1\mkern-2mu}\mkern 2mu}

\def\to{\rightarrow}

\def\0{\mymathbb{0}}
\def\1{\mathds{1}}
\def\z{\mymathbb{z}}
\def\a{\mymathbb{a}}
\def\b{\mymathbb{b}}
\def\c{\mymathbb{c}}
\def\d{\mymathbb{d}}
\def\x{\mymathbb{x}}
\def\y{\mymathbb{y}}

\def\N{\mathbb{N}}
\def\Z{\mathbb{Z}}
\def\Q{\mathbb{Q}}
\def\R{\mathbb{R}}
\def\C{\mathbb{C}}
\def\H{\mathbb{H}}

\def\D{\mathcal{D}}
\def\S{\mathcal{S}}
\def\E{\mathcal{E}}

\def\M{\mathcal{M}}
\def\s{\Sigma}
\def\r{\mathcal{R}}
\def\cpo{\mathbb{CP}^1}
\def\SUt{\mathrm{SU}(2)}

\def\SOtR{\mathrm{SO}_{3}\mathbb{R}}
\def\SLtZ{\mathrm{SL}_{2}\mathbb{Z}}
\def\SLtZn{\mathrm{SL}_{2}(\mathbb{Z}/n\mathbb{Z})}
\def\SLtR{\mathrm{SL}_{2}\mathbb{R}}
\def\SLtC{\mathrm{SL}_{2}\mathbb{C}}
\def\SLnZ{\mathrm{SL}_{n}\mathbb{Z}}

\def\SptnZ{\mathrm{PSp}_{2n}\mathbb{Z}}

\def\SptnZ{\mathrm{Sp}_{2n}\mathbb{Z}}

\def\GLtC{\mathrm{GL}_{2}\mathbb{C}}

\def\PSLtR{\mathrm{PSL}_{2}\mathbb{R}}
\def\PSLtZ{\mathrm{PSL}_{2}\mathbb{Z}}

\def\PGLtC{\mathrm{PGL}_{2}\mathbb{C}}

\def\slc{\mathfrak{sl}_2\mathbb{C}}

\def\glc{\mathfrak{gl}_2\mathbb{C}}
\def\glnC{\mathfrak{gl}_n\mathbb{C}}
\def\glnR{\mathfrak{gl}_n\mathbb{R}}

\def\mero{\mathfrak{M}}
\def\meros{\mathfrak{M}(\Sigma)}
\def\meroM{\mathfrak{M}(M)}

\def\meroA{\mathfrak{M}(M,\mathbb{A})}
\def\merofields{\mathfrak{M}\mathfrak{X}}

\def\merofieldsinvar{\mathfrak{M}\mathfrak{X}^\Gamma}

\def\meroforms{\mathfrak{M}\Omega}
\def\meroformss{\mathfrak{M}\Omega_\Sigma}

\def\invarmeroforms{\mathfrak{M}\Omega^\Gamma}

\def\invar{\mathfrak{M}^\Gamma}

\def\equi{\mathfrak{M}_\Gamma}
\def\equis{\mathfrak{M}_\Gamma(\Sigma)}
\def\equiA{\mathfrak{M}_\Gamma(M,\mathbb{A})}

\def\equicpo{\mathfrak{M}_\Gamma(\mathbb{CP}^1)}

\def\mop{\mathfrak{MD}}

\def\op{\mathfrak{RD}}
\def\opnc{\overbar{\mathfrak{RD}}}

\def\opeq{\mathfrak{RD}_\text{eq}}
\def\opnceq{\overbar{\mathfrak{RD}}_\text{eq}}

\def\opinv{\mathfrak{RD}^\text{inv}}
\def\opncinv{\overbar{\mathfrak{RD}}^\text{inv}}

\def\X{\langle X \rangle}
\def\Xinv{\langle X\rangle^\text{inv}}
\def\Xeq{\langle X\rangle_\text{eq}}

\def\killing{B}
\def\aut{\mathbb{A}\rm{ut}}
\def\auts{\mathbb{A}\rm{ut}\Sigma}
\def\autM{\mathbb{A}\rm{ut}\emph{M}}

\def\grassnR{\mathrm{Gr}_{n}\mathbb{R}^{2n}}

\def\siegeln{\mathbb{SH}_{n}}

\def\A{\mathbb{A}}
\def\units{\mathbb{A}^{\!\times}}
\def\apo{\mathbb{AP}^1}
\def\GLtA{\mathrm{GL}_{2}\mathbb{A}}
\def\PGLtA{\mathrm{PGL}_{2}\mathbb{A}}
\def\gla{\mathfrak{gl}_2\mathbb{A}}
\def\Id{\mymathbb{Id}}
\def\half{\textstyle\frac{1}{2}}
\def\Ad{\textnormal{Ad}}

\begin{document}

\title{Equivariant functions and rational differential operators}

\author{Michael~Deutsch}

\address{Departamento de Matematica, Instituto de Matematica\\ Universidade Federal do Rio de Janeiro\\ Caixa Postal 68530, Rio de Janeiro RJ 21941-909, Brasil.} \email{mdeutsch@im.ufrj.br}

\begin{abstract} We use contact geometry to describe the monoid of projectively equivariant meromorphic differential operators on a complex curve, quantization of which generalizes known constructions of classical equivariants to non-commutative function algebras in several variables.
\end{abstract}

 %32A17 (Primary) 58D07, 53C56 (Secondary)


\maketitle

Many constructions in geometry involve the fixed point set of a group action on sections of a fiber bundle. The present note is motivated by a desire to describe this set for a particularly simple class of (trivial) algebra bundles with non-compact group actions, the lowest rank instance of which is at once highly classical yet, to our knowledge, only explicitly resolved in special cases: Let $\s$ be a Riemann surface with function field $\mero$ and $\Gamma$ a group acting both on $\s$ by automorphisms $\Gamma \subset \auts$ and on the plane $\C$ (or sphere $\cpo \simeq \C \cup \infty$) by M\"obius transformations via a non-trivial projective representation $\rho : \Gamma \to \PGLtC$.
\begin{problem}\label{p1}
Construct all meromorphic functions $f : \s \to \C$ such that $f \circ A = \rho(A) \circ f$ for all $A \in \Gamma$.
\end{problem}
\noindent Examples of such \emph{equivariant functions}, the space of which we denote $\equis$, or simply $\equi$, are quite ubiquitous, occurring among Hauptmoduln of genus zero Fuchsian groups \cite{FJN94}, hyperelliptic covers of $\cpo$ \cite{GS05}, and Gauss maps of minimal immersions $\s \to \R^3$ \cite{Xu95}. Rational examples date back at least to Klein's famous study of the icosahedron \cite{Kl77}, where $A_5$-equivariants were obtained, including the (unique up to M\"obius conjugation) lowest-degree non-trivial example $$K(z) = \frac{z^{11}+66z^6-11z}{-11z^{10}-66z^5+1},$$ via a transvection-type construction defined on the ring of polynomial invariants of any finite $\Gamma \subset \GLtC$. While that construction already proves sufficient to resolve Problem \ref{p1} in the rational category \cite{DM89}, a non-trivial example with $\rho(\Gamma)$ non-compact seems only to have been observed as such nearly a century later, in the work of Heins \cite{He67} on elliptic function theory: Associating to a lattice $\Lambda = \langle z_1, z_2\rangle \simeq \langle \tau, 1\rangle \subset \C$ with modular constant $\tau = z_1/z_2 \in \H^2$ its Weierstrass zeta and eta functions, $\zeta_\Lambda : \C \to \C$ and $\eta_\Lambda : \Lambda \to \C$,  $$\zeta_\Lambda(z) = \frac{1}{z} + \sum\limits_{w \in \Lambda^*}(\frac{1}{z-w} + \frac{1}{w} + \frac{z}{w^2}), \quad \eta_\Lambda(w) = \zeta_\Lambda(z+w) - \zeta_\Lambda(z),$$ Heins noticed that the function $H : \H^2 \to \C$ defined by $H(\tau) = \eta_\Lambda(\tau)/\eta_\Lambda(1)$ is equivariant with respect to the modular group $\PSLtZ$, and established mapping properties of general modular equivariants in order to study the pseudo-periods of $\zeta$. Brady \cite{Br71} generalized this observation and solved Problem \ref{p1} for an arbitrary Fuchsian group $\Gamma$ acting on $\H^2$ by determining an algebraic embedding of the field of $\Gamma$-automorphic functions $\invar \hookrightarrow \equi$, a result which Smart \cite{Sm72} generalized to the case of an arbitrary Kleinian group with domain of discontinuity $\s \subset \cpo$. The finite and modular cases have reappeared in multiple contexts, such as numerical analysis where Doyle and McMullen rediscovered $K$ and proved that its $20$ critical points form superattracting $2$-cycles, leading to a root-finding algorithm for the quintic \cite{DM89}, and in differential geometry and physics, where general rational equivariants correspond to singular points in the moduli space of monopoles on $\R^3$ \cite{Do84}, \cite{Ja00}, \cite{HMS98}; similarly, Nahm rediscovered $H$ in the context of quantum field theory \cite{Na85}, number theoretic aspects of which have been investigated by Sebbar \cite{SS12} and, indirectly, Duke \cite{Du05}.\\

%interpreting elements of $\equicpo$ as homogeneous forms on $\C^2$ invariant with respect to the double cover $\tilde{\Gamma} \subset \SLtC$,

%One such context where straight-forward extensions of classical methods (the Klein construction, Reyonlds operators, Molien series) have already been successfully implemented is that of morphisms \cite{CS16} between projective varieties (especially projective space itself \cite{dFH18}) equivariant with respect to compact (especially finite \cite{Pr12}) group actions, but suitable adaptations of this machinery do not seem to have been developed for function spaces with non-compact actions. We supply such an adaptation when the source and target, $\s$ and $\C$, are replaced with an arbitrary complex $\Gamma$-manifold and unital complex Banach algebra, $M$ and $\A$. A morphism $\rho : \Gamma \to \PGLtA$ determines an action on the algebra of $\A$-valued functions $\meroA$

While applications are evidently diverse, solutions to special cases of Problem \ref{p1} suggest a general principle, essentially due to Klein, whereby the task of finding equivariant objects should be translated, by some device, into that of finding suitable \emph{in}variant objects. Motivated by the extent to which the constructions of Heins and Brady might generalize to moduli of curves of higher genus and abelian varieties, we aim to apply this principle in a higher dimensional, non-commutative context: Given a complex $\Gamma$-manifold $M$ and unital complex Banach algebra $\A$, a morphism $\rho : \Gamma \to \PGLtA$ defines a $\Gamma$-action
\begin{align*}
\Gamma \times \meroA \; & \to \; \meroA \\
(A,f)\; & \mapsto \; \rho(A) \circ f \circ A^{-1}
\end{align*}
on the algebra $\meroA$ of $\A$-valued functions, the fixed point set of which we denote $\equiA$. Our main observation is a differential incarnation of Brady's embedding:

\begin{thm}\label{deform}
An equivariant function $f_0 \in \equiA$ such that $\dot{f}_0^{-1} = X(f_0)^{-1} \in \meroA$ for some $\Gamma$-invariant meromorphic vector field $X \in \merofieldsinvar$ admits an explicit deformation $\mathfrak{F} \subset \equiA$ containing a subfamily $\{f_t\}_{t \in \mathcal{U}} \subset \mathfrak{F}$ parameterized by a neighborhood of the zero function $0 \in \mathcal{U}$ in the field of scalar invariants $\mathcal{U} \subset \invar$.
\end{thm}

\noindent As we will see, in the classical case $M = \s$ with $\A = \C$, the neighborhood $\mathcal{U}$ extends to $\invar$ and reduces Problem \ref{p1} to that of finding just one $\Gamma$-equivariant function. That such a deformation should exist is intuitive from the point of view of the equivariant Reynolds operator (the use of which trivially resolves the matter for $\Gamma$ compact \cite{dFH18}), but apart from that method, all those we are aware of can be obtained from this deformation. Our approach to establishing it is inspired on the one hand by the rather enviable monoidal structure in the rational case, whereby the action of $\mero_{\rho(\Gamma)}(\cpo)$ on $\mero_\Gamma(\s)$ produces new equivariants for $\rho(\Gamma)$ compact, and on the other by \cite{Ra56}, where Rankin determined which polynomial expressions in the derivatives of an automorphic form are again automorphic. We will determine which rational expressions in derivatives of an equivariant function are again equivariant, an analogue of Problem \ref{p1} at the level of differential operators, where we call a quotient $\r = \mathcal{P}/\mathcal{Q} \in \op$ of meromorphic operators $\mathcal{P}$ and $\mathcal{Q}$ a \emph{rational operator}:
\begin{problem}\label{p2}
Construct all rational operators $\r : \mero \to \mero$ such that $\r \circ T = T \circ \r$ for all $T \in \PGLtC$.
\end{problem}
\noindent The monoid $\opeq$ of these \emph{projectively equivariant} operators is the natural counterpart of the algebra $\opinv$ of invariant operators, and we show that the two are indeed related according to Klein's principle, using the contact geometry of the M\"obius group as the translation device. Briefly, in sections \ref{SectionCurve} and \ref{SectionOperator}, we show how a certain duality notion \cite{KUY03} for Legendrian curves in $\PGLtC$ induces an involution on an open subset of pairs $\mero \times \meroforms$ of meromorphic functions and $1$-forms, the first component $\D$ of which is a manifestly equivariant analogue of the second $\S$, the Schwarzian derivative. The key observation is that the former admits a trivial holomorphic deformation with base $\opinv$, which exhausts $\opeq$ over a curve $\s$ and solves Problem \ref{p2}. In section \ref{SectionClassical} we establish the relation between the resulting solution to Problem \ref{p1} and the classical methods, and in section \ref{SectionNC} define a larger algebra $\opnc$ of \emph{non}-commutative rational operators on an arbitrary manifold $M$. We then adapt the procedure above to the generalized M\"obius group $\PGLtA$ which, over a curve, associates to each $\PGLtC$-equivariant rational operator $\r \in \opeq$ a unique $\PGLtA$-equivariant ``quantization'' $\overbar{\r} \in \opnceq$ -- an \emph{a priori} non-commutative version which reduces to $\r$ in the ``classical limit'' $\opnc \to \op$.

\begin{thm}\label{oper1}
Given a meromorphic vector field $X$ on $M$, let $\X \subset \opnc$ denote the free associative skew-field over $\mero$ generated by iterates of $X$, and $\Xeq = \X \cap \opnceq$ the submonoid of $\PGLtA$-equivariant elements. Then $$\E\, \mapsto \, \Id + X\,[\E\circ_q (X\circ\varphi_X + \varphi_X^2)+\varphi_X]^{-1}$$ maps $\X$ onto $\Xeq - \{\Id\}$, where $\varphi_X = -\half X^{-1}(X\circ X)$ and $\circ_q$ is defined inductively by the relation $X \circ_q \mathcal{H} = X \circ \mathcal{H} - [\varphi_X,\mathcal{H}]$, with $[\cdot,\cdot]$ denoting the multiplicative commutator.
\end{thm}

\noindent As we will make explicit, this represents a sort of Cayley transform which again embeds a certain space of invariant objects into that of equivariant ones. With this in hand, if $M$ admits an invariant meromorphic vector field $X$, a subfamily of operators stabilizing $\equiA$ can be extracted from $\Xeq$, the action of which establishes Theorem \ref{deform} and produces a large class of examples with $\rho(\Gamma)$ non-compact. While the strictly function-theoretic content is all verifiable directly without appeal to the contact structure on $\PGLtC$, we hope the reader will appreciate the efficiency this geometric link affords.

\section{Legendrian curves in $\PGLtC$}\label{SectionCurve}

Consider the following geometric attempt to build a M\"obius equivariant sheaf morphism $\mero \to \mero$ on the function field of a Riemann surface $\s$: Given a non-constant function $f \in \mero$, regard its values as distributed in the ideal boundary $\partial_\infty\H^3 = S^2$ of the Poincaré ball $\H^3 \subset \R^3$ with hyperbolic metric $g$. Let $U \subset \s$ be an open set on which $f|_U$ is univalent, and choose an oriented immersion $x = x_f : U \subset \s \to \H^3$ with (positive) \emph{hyperbolic Gauss map} $f$, that is, for any $p\in U$ the oriented hyperbolic geodesic $\gamma_{x(p),N}$ normal to $x(U)$ at $x(p)$ meets the boundary at $f(p) = \lim_{t \to \infty}\gamma_{x(p),N}(t)$. Then the reflection $$f(p) \mapsto \hat{f}(p) := \lim_{t \to -\infty}\gamma_{x(p),N}(t)$$ of $f(U) \subset S^2$ through $x(U)$ defines a new function $\hat{f} : U \to S^2$, the \emph{negative} hyperbolic Gauss map. Since ambient-congruent immersions have M\"obius-equivalent Gauss maps, if this choice could be made compatibly over \emph{all} such function elements $(f,U)$, then $f \to x_f$ would behave as a morphism of sheaves with $x_{T\circ f} = T \circ x_f$ for all $T \in \PGLtC \simeq \textnormal{Isom}(\H^3)$, and the induced morphism $f \mapsto \hat{f}$ would be intrinsically $\PGLtC$-equivariant.\\

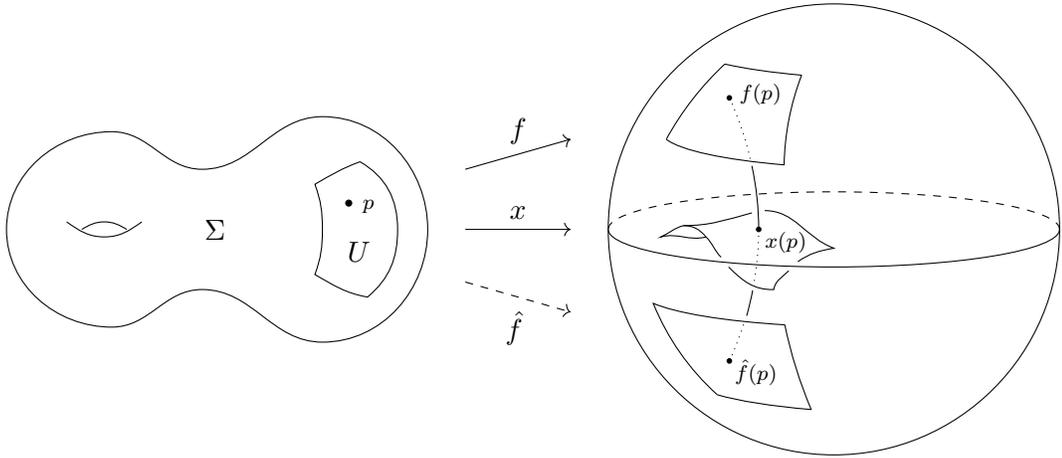
\begin{tikzpicture}
%surface \Sigma
\draw (5.6,0) to [out=90,in=0] (4.2,1.5) to [out=180,in=0] (2.6,0.8) to [out=180,in=0] (1.4,1.3) to [out=180,in=90] (0,0) to [out=-90,in=180] (1.4,-1.3) to [out=0,in=180] (2.6,-0.8) to [out=0,in=180] (4.2,-1.5) to [out=0,in=-90] (5.6,0);
\draw (0.8,0.1) to [out=-35,in=180] (1.3,-0.1) to [out=0,in=215] (1.8,0.1);
\draw (1,0) to [out=30,in=180] (1.3,0.1) to [out=0,in=150] (1.6,0);
\node[right] at (2.5,0) {$\Sigma$};

%open set U
\draw (4.1,0.6) to [out=30,in=200] (4.7,0.9) to [out=-30,in=90] (5.2,0) to [out=-90,in=40] (4.8,-0.9) to [out=170,in=-30] (4.1,-0.6) to [out=70,in=-90] (4.2,0) to [out=90,in=-70] (4.1,0.6);
\draw[fill] (4.55,0.35) circle [radius=0.04]; \node[right] at (4.6,0.3) {${\scriptstyle p}$}; \node[right] at (4.4,-0.3) {$U$};

%maps f, x, \hat{f}
\draw[->] (6.1,0) -- (7.5,0); \node[above] at (6.8,0) {$x$}; \draw[->] (6.1,0.8) -- (7.5,1.2); \node[above] at (6.8,1) {$f$}; \draw[dashed, ->] (6.1,-0.7) -- (7.5,-1.1); \node[below] at (6.75,-0.95) {$\hat{f}$};

%sphere
\draw (11,0) circle (3cm);
%\equator
\draw (8,0) arc (180:360:3cm and 0.5cm); \draw[dashed] (14,0) arc (0:180:3cm and 0.5cm);

%geodesic \gamma
\draw[dotted] (9.61,1.75) arc (25:12.5:4cm); \draw (10,0) arc (0:12.2:4cm);
\draw[dotted] (10,0) arc (0:-10:4cm); \draw[dotted] (9.61,-1.75) arc (-25:-17:4cm);
\draw (9.84,-1.15) arc (-15.8:-10.8:4cm);
%f(p)
\draw[fill] (9.61,1.75) circle [radius=0.03]; \node[right] at (9.61,1.8) {${\scriptstyle f(p)}$};
%x(p)
\draw[fill] (10,0) circle [radius=0.03]; \node[below right] at (9.95,0.1) {${\scriptstyle x(p)}$};
%\hat{f}(p)
\draw[fill] (9.61,-1.75) circle [radius=0.03]; \node[right] at (9.55,-1.9) {${\scriptstyle\hat{f}(p)}$};

%f(U)
\draw (9.55,2.2) arc (133:151.2:4cm); \draw (9.55,2.2) arc (216:242.1:3cm and 0.5cm);
\draw (8.77,1.2) arc (192:243:3cm and 0.5cm); \draw (10.57,2.05) arc (159:179:3.5cm);

%x(U)
\draw (10.52,-0.43) to [out=35,in=190] (11,-0.25) to [out=160,in=0] (10.08,0.25); \draw (9.91,0.22) to [out=195,in=0] (9,0.05) to [out=195,in=25] (8.7,-0.1) to [out=15,in=195] (9.3,-0.1); \draw (9,0.05) to [out=0,in=130] (9.3,-0.1) to [out=-45,in=135] (9.5,-0.37);  \draw (9.6,-0.5) to [out=-45,in=182] (10.2,-0.8) to [out=80,in=225]  (10.4,-0.55);

%hat{f}(U)
\draw (9.45,-2.2) arc (228:202:3.3cm); \draw (10.7,-2.4) arc (203:191.7:6cm); \draw (8.59,-0.98) arc (201:249.7:3cm and 0.5cm); \draw (9.45,-2.2) arc (210:243.2:3cm and 0.5cm);
\end{tikzpicture}

\vspace{0.5cm}
\noindent Promoting this still ambiguous, local concept to a well-defined global operator $\mero \to \mero$ leads immediately to the complex geometry of $\PGLtC$ itself by two classical theorems:\\

\noindent Bianchi \cite{Bi24}: \emph{Reflection $f(U) \subset \hat{\C} \to \hat{f}(U) \subset \hat{\C}$ through $x_f(U)$ is holomorphic if and only if the induced metric $x_{f}^*g$ is intrinsically flat.}\\ %that is, the Gaussian curvature of the induced metric $x_f^*g$ vanishes identically,

\noindent G{\'a}lvez, Mart{\'{\i}}nez, Mil{\'a}n \cite{GMM00}: \emph{A front $x : \s \to \H^3$ is flat if and only if, locally and away from singular points, $x$ lifts to a Legendrian immersion $L : U \subset \s \to \PGLtC$ holomorphic with respect to the conformal structure determined by the second fundamental form of $x$.}\\

The M\"obius group $\PGLtC = \{[T] \in \mathbb{P}(\glc)  \, | \, \det T \ne 0 \}$ and its double cover $\SLtC = \{T \in \glc \, | \, \det T = 1\}$ share the Lie algebra $\slc$ of trace-free matrices. Regarded as a left-invariant vector field, the element $R = \textnormal{diag}(1,-1) \in \slc$ generates a holomorphic flow, right-multiplication by the integral curve $\C^* \simeq \{[\textnormal{diag}(a,a^{-1})]\, | \, a \ne 0 \}$, and the quotient $\PGLtC/\C^*$ is naturally identified with the space of point pairs in the Riemann sphere $\cpo$,
\begin{align*}
 \pi : \PGLtC &\to \cpo \times \cpo - \Delta\\
  [T] &\mapsto ([T_1],[T_2])
\end{align*}
where $[T_1]$ and $[T_2]$ denote the column spaces of either representative matrix $T \in \SLtC$ of $[T] \in \PGLtC$. Consider the holomorphic distribution of tangent complex 2-planes $p \mapsto R_p^\bot \subset T_p\PGLtC$ given by the holomorphic Killing-orthogonal complements to the fibers of $\pi$. Evidently $R_p^\bot = \textnormal{ker}\,\xi_p$, where $\xi$ is the 1-form appearing on the diagonal of the left-invariant Maurer-Cartan form $$\omega = \begin{pmatrix} \xi & \hat{\theta} \\ \theta & -\xi \end{pmatrix} : T\PGLtC \to \slc,$$ whose entries provide a basis for the dual Lie algebra. The structure equation $d\omega + \omega \wedge \omega = 0$ implies that $\xi \wedge d\xi = \xi \wedge \theta \wedge \hat{\theta}$ is a holomorphic volume form, so this distribution is everywhere non-integrable in the sense of Frobenius, i.e. a holomorphic \emph{contact structure}. In addition to left-invariance, the distribution is manifestly \emph{right}-invariant by the Reeb flow $\C^* \simeq \textnormal{exp}([R])$ and by the involution $$N = \begin{bmatrix} 0 & 1 \\ 1 & 0 \end{bmatrix}$$ interchanging the factors of the projection $\pi : \PGLtC \to \cpo \times \cpo - \Delta$, all of which are holomorphic isometries of the (re-scaled) Killing form $\killing = \xi^2 + \theta\hat{\theta},$ viewed as a complex Riemannian metric of constant curvature.\\

By a \emph{holomorphic curve} in the M\"obius group we mean any non-constant holomorphic map from a Riemann surface $L : \s \to \PGLtC$ defined on the complement of a discrete subset of ends. Such a map yields a pair of meromorphic functions $f$ and $\hat{f}$, the components of the projected curve $\pi(L) = (f,\hat{f})$ which extends holomorphically to the diagonal $\Delta \subset \cpo \times \cpo$ across appropriate ends of $L$. A curve $L$ is called \emph{Legendrian}, or \emph{contact}, if it is tangent to the contact distribution at all regular points, $L^*\xi = 0$. In this case, the remaining entries $\theta$ and $\hat{\theta}$ of the pullback $L^*\omega$ of the Maurer-Cartan form are called the \emph{left} and \emph{right} \emph{canonical 1-forms of L}, at least one of which must be non-zero, which locally determine $L$ up to constant left $\PGLtC$-multiple. Lifting $L$ to $\SLtC$ locally, the tangent vector can be viewed as a matrix-valued map $dL : U \subset \s \to \glc$, the column spaces of which coincide with those of $L$ by the contact condition, so $\pi(L)$ may be thought of as a kind of Gauss map: We call $f$ and $\hat{f}$ the \emph{left} and \emph{right} \emph{Gauss maps} of $L$, at least one of which must be non-constant, and which determine $L$ locally up to a constant right $\C^*$-multiple \cite{KUY03}. The data pair $(f,\theta)$ however determines the curve globally and uniquely:

\begin{lem}\label{lift}
Given $f \in \meros$ non-constant and $\theta \in \meroformss$ non-zero, there exists a unique Legendrian curve $L_{f,\theta} : \s \to \PGLtC$ with left Gauss map $f$ and left canonical form $\theta$.

\begin{proof} Denoting $\dot{f} = df/\theta$ and defining auxiliary meromorphic functions $\varphi = -\half (\ddot{f}/\dot{f})$ and $Q = \dot{\varphi}+\varphi^2$, the first-order linear system $dL = L\omega$ can locally be re-written as a single second-order condition: Any point in the complement $\s_0 = \s - \textnormal{Supp}(\dot{f})$ of the divisor $(\dot{f})$ has a neighborhood $U \subset \s_0$ where $f = [\psi]$ lifts to a solution $\psi : U \to \C^2$ of the complex Schr\"odinger equation $\ddot{\psi}+Q\psi = 0$. Since such lifts agree up to a constant $\C^*$-multiple, the projection $L_{f,\theta} = [\psi,\dot{\psi}]$ of the matrix function with columns $\psi$ and $\dot{\psi}$ extends to a globally defined meromorphic curve $L_{f,\theta} : \s \to \PGLtC$, easily seen to be Legendrian.
\end{proof} %(namely $\psi = (f\dot{f}^{-1/2}, \dot{f}^{-1/2})^t$ where $\dot{f}^{1/2} : U \to \C^*$ is any branch of the square-root of $\dot{f}$ on $U$ simply connected)
\end{lem}

As a consequence, the space $\mathfrak{L}_\s$ of Legendrian curves parameterized by $\s$ can be identified with a subset of $\meros \times \meroformss$. When $\theta$ is fixed, this gives a lifting procedure whereby a non-constant $f : \s \to \cpo$ has a unique horizontal lift $L : \s \to \PGLtC$ relative to $\theta$,

\begin{center}\begin{tikzcd}[row sep=large]
& \PGLtC \arrow{dl}[swap]{\pi_1}\arrow{dr}{\pi_2} & \\
\cpo & \s \arrow{l}[swap]{f}\arrow[u, dashed, "L_{f,\theta}" description]\arrow[r, dotted, "\hat{f}"] & \cpo
\end{tikzcd}\end{center}

\noindent  where $\pi_1$, $\pi_2$ is the double fibration of $\PGLtC$ by null planes\footnote{That is, totally geodesic complex hypersurfaces on which the Killing form is everywhere degenerate.} intersecting in the geodesic fibers of $\pi = (\pi_1,\pi_2)$. The subspace of curves $L_{f,\theta}$ (and $L_{f,\theta}N$) contained entirely within complex geodesics corresponds to the subspace of data pairs satisfying a non-linear first order equation, $\mathfrak{N}_\s = \{\, (f,\theta)\, | \, \dot{f} = (a_1f+a_2)^2, \; (a_1, a_2) \in \C^2-\{0\} \},$ which are precisely the Legendrian curves with $\hat{f}$ (resp. $f$) constant.

\section{The generating operators $\D$ and $\S$}\label{SectionOperator}

Let $X \in \merofields$ be a meromorphic vector field on $\s$, $\theta \in \meroforms$ its dual 1-form $\theta(X) = 1$ and $\dot{f} = X(f)$, and denote the $n$-fold product $X^n : f \mapsto \dot{f}^n$ and operator composition $X_n : f \mapsto f^{(n)}$, $$X_n = \underbrace{X \circ X \circ \ldots \circ X}_{n \text{ times}}.$$ The algebra of \emph{meromorphic differential operators} on $\s$ is the infinite polynomial ring\footnote{Our terminology differs from that of Donin and Khesin \cite{DK07} who use this phrase to refer to the subset $MDO = \{D = \sum_{i=1}^{n} a_iX_i \;|\; a_i \in \mero, n\in\N\}$ with Lie bracket $[D_1,D_2] = D_1\circ D_2 - D_2\circ D_1$.} $$\mop = \mero[X_0,X_1,X_2,...]$$ over the function field $\mero$, where $X_0 = \Id_{\mero}$.  As a composition ring, $(\mop,+,\cdot,\circ)$ is an integral domain, and its definition is independent of the choice of $X$. We call its field of fractions the algebra of \emph{rational differential operators} on $\s$, $$\op = \mero(X_0,X_1,X_2,...).$$ This larger algebra is still a composition ring, an element $\r \in \op$ acting on the open set of $\mero$ where its denominator does not vanish identically, thus defining a left monoid pseudo-action of $\op$ on $\mero$ which extends both the (genuine) submonoid action by $\mop \subset \op$ and that of its maximal subgroup $\PGLtC \subset \op$. Pre-composition also defines a right action on the operators $\op \otimes E$ taking values in the space of meromorphic sections of a given holomorphic vector bundle $E \to \s$.\\

Consider the pair of operators $\D_\theta \in \op$ and $\S_\theta \in \op \otimes \Omega$ (dually, $\D_X$ and $\S_X$) with respect to $\theta \in \meroforms$ ($X \in \merofields$) defined by the components $$\D_\theta f = \pi_2 \circ L_{f,\theta} = \hat{f}, \quad \quad \S_\theta f = L_{f,\theta}^{\,*} \, \killing(X,\cdot) = \hat{\theta}$$ of the transform $L \mapsto LN$ on $\mathfrak{L}_\s$ mapping the left data set to the right.

\begin{prop}\label{operators} $\D$ is projectively equivariant, $\mathcal{S}$ projectively invariant, and for $(f,\theta) \notin \mathfrak{N}_\s$ $$\D_{\S_\theta f} \circ \D_\theta (f) = f, \quad \quad \S_{\S_\theta f} \circ \D_\theta (f) = \theta.$$
\begin{proof}
Given $T \in \PGLtC$, the proof of Proposition \ref{lift} shows that $L_{T\circ f,\theta} = TL_{f,\theta}$. Then $\pi_2 \circ T = T \circ \pi_2$ implies $\D \circ T = T \circ \D$, while left-invariance of the Maurer-Cartan form implies $\S \circ T = \S$. The identities relating $\D$ and $\S$ are equivalent to the assertion that $(f,\theta) \mapsto (\hat{f},\hat{\theta})$ is an involution on the set of non-geodesic Legendrian curves.
\end{proof}
\end{prop}

Equivariance $L_{T\circ f,\theta} = TL_{f,\theta}$ also holds on the degenerate set $\mathfrak{N}_\s$, and for any $\theta \in \meroforms$,
$$\D_{\theta}f = T \circ \D_{\theta}g \;\;\; \Leftrightarrow \;\;\; \S_{\theta} f = \S_{\theta} g \;\;\; \Leftrightarrow \;\;\; f = T \circ g \;\;\; \Leftrightarrow \;\;\; \S_{dg} f \equiv 0 \;\;\; \Leftrightarrow \;\;\; \D_{dg}f \equiv \textnormal{constant}.\footnote{It is convenient here to regard the point at infinity $\infty$ as a constant function.}$$

\noindent Thus $\S$ is a $1$-form-valued incarnation of the classical Schwarzian derivative, $\D$ its natural equivariant counterpart. The former can be defined more invariantly as the induced metric $L_{f,\theta}^{\,*} \, \killing \in \op \otimes \Omega^2$, usually interpreted as a meromorphic quadratic differential, or more traditionally as an operator $S \in \opinv$, typically defined with respect to a function $z : U \subset \s \to \C$ as the coefficient $S_z f\, dz = \S_{dz}f$, although we will abuse the distinction and use the same notation for both in the sequel, as well as $D_z f = \D_{dz} f$. The latter is in fact a global instance of our earlier local attempt to design manifestly equivariant morphisms:

\begin{example} Up to a normal shift \cite{GMM00}, an oriented flat immersion $x$ is determined by its hyperbolic Gauss map $f$ and \emph{Hopf differential} $\mathcal{Q}$, the $(2,0)$-part of the induced metric with respect to the conformal structure determined by the second fundamental form of $x = x_{f,\mathcal{Q}}$. Since M\"obius transformations of the boundary extend to isometries of the interior, the map $f \to x_{f,\mathcal{Q}_f}$ on germs will be projectively equivariant as soon as $f \to \mathcal{Q}_f$ is projectively invariant. Given any fixed 1-form $\theta$, the simplest choice $\mathcal{Q}_f = S_\theta f\,\theta^2$ then makes reflection $\mero \to \mero$ of $f$ through $x_{f,\mathcal{Q}_f}$ not only equivariant, but globally defined, coinciding precisely with the $\D$ operator, $f \mapsto \D_\theta f$.\end{example}

Local mapping properties of $\S$ are well-known from classical complex analysis, while those of $\D$ can be derived from the \emph{pre-Schwarzian} $\varphi_z f = -\half (\ddot{f}/\dot{f})$ by writing  $$D_z f = f + {\textstyle\frac{df}{\varphi_z f\,dz}}, \quad S_z f = {\textstyle\frac{d}{dz}}(\varphi_z f) + (\varphi_z f)^2.$$ The change of variable formulae for $D_z$ and $S_z$ then follow from $$(\varphi_w f)dw = (\varphi_z f - \varphi_z w)dz,$$ allowing $\varphi_zdz$ and $S_zdz^2$ to be interpreted as cocycles on the pseudo-Lie group of local coordinate transformations on $\s$ \cite{Gu67}. One also infers the behavior of $D_z f$ and $S_z f$ at ramification points of $f$:

\begin{prop} On a coordinate neighborhood $(U,z)$, the relations $\textnormal{Ram}(D_z f) = \textnormal{Ram}(f)$ and $\textnormal{Pol}(S_z f) = -2\,\textnormal{Supp}\,\textnormal{Ram}(f)$ hold between the ramification and polar divisors. \end{prop}

Similarly, global mapping properties of $\S$ are well-known from classical Teichmüller theory, while those of $\D$ can be derived from the duality described in Proposition \ref{operators}, such as the surjectivity of the map
\begin{align*}
\mero \times \meroforms & \; \to \; \mero\\
(f,\theta) & \; \mapsto \; \D_\theta f,
\end{align*}
under which the pre-image of a point $\hat{f} \in \mero$ is $\{(\D_\vartheta \hat{f}, \S_\vartheta \hat{f}) \, | \, \vartheta \in \meroforms \}.$
However, the topology of $\s$ may obstruct surjectivity on restrictions of the form
$$\mero \times \{\theta\} \to \mero,\quad \{f\} \times \meroforms \to \mero,$$
since given a fixed 1-form $\theta$, a $\D_\theta$-primitive of $\hat{f} \in \mero$ is globally well-defined if the differential equation $\D_{\theta} f = \hat{f}$ in $f$ has trivial monodromy on $\s$, while a given $f$ is a $\D_\theta$-primitive of $\hat{f}$ for some $\theta$ if the period map from the fundamental group $\pi_1(\s) \to \C$ $$[\gamma]\; \mapsto \;\frac{1}{\pi i}\oint_{\gamma}\frac{df}{f-\hat{f}}$$ takes values in $\Z$ (in which case $\theta = \textnormal{exp}[2\int(f-\hat{f})^{-1}df]df$ is well-defined), but generally not otherwise.\\

Thus to implement $\D$ in resolving Problem \ref{p1} on surfaces of arbitrary topological type, it is desirable to solve Problem \ref{p2} in the same generality. Note that the analogue of the latter problem for invariant operators is trivially resolved by Cartan's theory of moving frames \cite{Ol18}: The uniqueness of $L_{f,\theta}$ allows it to be interpreted as the end result of a reduction of $\PGLtC$-frames over $f$ with respect to the Borel subgroup stabilizing $[1,0]^t \in \cpo$. Since $S$ is the unique non-constant entry of $(L_{f,\theta}^{\,*}\omega)/\theta$, it not only represents the lowest-order non-constant element of $\opinv$, but in fact the ``complete set'' of projective invariants in dimension 1, determining $f$ up to M\"obius transformation. Consequently, any invariant of the prolonged action of $\PGLtC$ on the holomorphic jet bundle of $\s$ is a rational function of its derivatives:

\begin{prop}\label{invop} Derivatives of the Schwarzian generate the subalgebra of $\PGLtC$-invariant rational operators, in the sense that given any non-zero meromorphic vector field $X \in \merofields$, pre-composition $\r \mapsto \r \circ S_X$ maps $\op$ to $\opinv$ surjectively.
\end{prop}

\noindent While Cartan's theory does not directly apply to equivariant operators, an indirect application similarly singles out the $\D$ operator as the lowest-order non-trivial element of $\opeq$. Though lacking the algebra structure of $\opinv$, $\opeq$ does possess a monoidal structure, of which iterates of $\D$ determine a discrete submonoid. The structure of the latter begins to appear already at the second iteration, $$\D_X \circ \D_X = \Id + X/(\mathcal{H}+\varphi_X),$$ which resembles the change of variable formula for $\D$ but with $\mathcal{H} = (X\circ S^{-1})^{-1} \in \opinv.$ Similarly, it is not difficult to verify that the deformation $\gamma(t) = \Id + X/(t\mathcal{H}+\varphi_X)$ of $\gamma(0) = \D$ into $\gamma(1) = \D \circ \D$ remains entirely in $\opeq$ for $0 \le t \le 1$, from which one might suspect that all equivariant operators could be so expressed for some $\mathcal{H} \in \opinv$. The next result establishes that suspicion affirmatively which, combined with Proposition \ref{invop}, explicitly resolves Problem \ref{p2} and identifies $\opeq$ with an infinite-dimensional vector space plus a ``point at infinity'':

\begin{thm}\label{operator}
Deformations of the $\D$ operator generate the submonoid of $\PGLtC$-equivariant rational operators, in the sense that given any non-zero meromorphic vector field $X \in \merofields$, the deformation $\Phi_X : \op \to \op$ of $\D$ defined by $$\mathcal{H} \, \mapsto \, \Phi_X\mathcal{H} = \Id + X/(\mathcal{H}+\varphi_X)$$ maps $\opinv$ to $\opeq - \{\Id\}$ bijectively.

\begin{proof}
Since $\Phi_X\mathcal{H}(f) = L_{f,\theta} \circ \mathcal{H}(f)$ with $\theta(X) = 1$ and $L_{f,\theta}$ acting point-wise, the equivariance property of $L$ implies that $\Phi_X(\opinv) \subset \opeq$, and by projective invariance of the cross-ratio $[\cdot, \cdot, \cdot, \cdot]$, the map $\E \mapsto [\Phi(\infty),\Phi(0),\Phi(1),\E]$ inverts $\Phi$, where $\Phi(\infty) = \Id$.
\end{proof}
\end{thm}

A small subfamily in $\opeq$ is now sufficient to solve Problem \ref{p1}: Identifying an element $h \in \mero$ with the operator $\r : f \mapsto h$, the function field itself is viewed as a subspace $\mero \subset \opinv$, that of the constant operators. We have $(f,\theta) \in \equi \times \invarmeroforms$ if and only if the curve $L_{f,\theta} : \s \to \PGLtC$ is left $\rho(\Gamma)$-equivariant ($L \circ A = \rho(A)L$ for all $A \in \Gamma$) if and only if $h \, \mapsto \, \Phi_Xh(f)$ maps $\invar$ to $\equi - \{f\}$ with inverse $g \mapsto L_{f,\theta}^{-1}\circ g$.

\begin{cor}\label{bijection}
Given $f \in \equi$ and $g \in \invar - \{0\}$, the map $$h \, \mapsto \, f + \dot{f}\,(h - \half \dot{f}^{-1}\ddot{f})^{-1}$$ defines a bijection between $\invar$ and the punctured space $\equi - \{f\}$, where $\dot{f} = df/dg.$
\end{cor}

\section{Classical equivariant functions: Single variable commutative theory}\label{SectionClassical}

In light of Corollary \ref{bijection}, it is reasonable to expect that the classical methods alluded to in the introduction should appear as special cases of the $\D$ operator and its deformations. We now demonstrate precisely this in a series of examples, the collection of which exhausts (our knowledge of) the history of Problem \ref{p1} and suggests a unified perspective on the subject, surveying some properties of equivariants themselves in the process. Specifically, we will derive from $\D$ a family of operators on the ring of invariants, a special case of which we call the $\phi$-\emph{operator}, of interest in its own right, to which the classical constructions can be reduced.\\

This setting concerns surfaces conformal to a subset of the sphere $\s \subset \cpo$ acted on by a group of M\"obius transformations $\Gamma \subset \GLtC$. A (meromorphic) $\Gamma$-\emph{automorphic form}, or simply \emph{invariant}, of \emph{weight} $k$ and \emph{character} $\chi : \Gamma \to U(1) \subset \C$ is then a meromorphic function $\alpha : \s \to \C$ such that $$\alpha(A \circ z) = \chi(A)(cz+d)^k \alpha(z) \;\;\; \textnormal{for all} \;\;\; A = \begin{pmatrix} a & b \\ c & d \end{pmatrix} \in \Gamma,$$ and an \emph{automorphic function} or \emph{absolute invariant} a weight zero form with trivial character. Suppressing both the group and the character in the notation, the vector space of weight $k$ automorphic forms is denoted $\M_k$, and the direct sum $\M_* = \bigoplus_{k \in \Z}\M_k$ is a graded ring with multiplication $[\alpha,\beta]_0 = \alpha \beta$ when $\chi = 1$. The derivative $\alpha \mapsto \alpha' = \textstyle\frac{d\alpha}{dz}$ does not preserve $\M_*$, but for $\alpha \in \M_k$, $\beta \in \M_l$, the bracket $[\alpha,\beta]_1 = k\alpha \beta' - l\alpha'\beta$ does, with respect to which $\M_*$ forms a graded Poisson algebra.  These two products constitute the zero- and first-order \emph{Rankin-Cohen brackets} \cite{Ra56}, \cite{Co75}, $$[\alpha,\beta]_n = \sum\limits_{m=0}^{n} (-1)^m\binom{n+k-1}{n-m}\binom{n+l-1}{m}\alpha^{(m)}\beta^{(n-m)}$$
defining $(-1)^n$-symmetric bi-differential operators $$[\cdot,\cdot]_n : \M_k \times \M_l \to \M_{k+l+2n},$$ providing coefficients \cite{Za94} of a genuine deformation quantization $\alpha *_{\hbar} \beta = \sum_{n=0}^\infty [\alpha,\beta]_n \hbar^n$ of $(\M_*,[\cdot,\cdot]_1)$. For a finite group $\Gamma \subset \GLtC$ acting on $\s = \hat{\C}$, these are essentially the classical \emph{transvectants} \cite{Gor87} produced by the Cayley $\Omega$-process on the subring $\mathcal{P}_* = \C[z] \cap \M_*$ of polynomial invariants. Analogously, for a finite-index congruence subgroup $\Gamma \subset \SLtZ$, a \emph{modular form} is a holomorphic $\Gamma$-automorphic form on $\s = \H^2$ extending to the cusps $C \subset \Q \cup \{\infty\}$ of $\Gamma$, so that a weight $2k$ element $\alpha$ corresponds to a $k$-\emph{differential} $\alpha(z)dz^k$, that is, a holomorphic section of the $k^{th}$ power of the canonical bundle $\Omega_{\overline{\H^2/\Gamma}}$ over the compactification of the quotient $\H^2/\Gamma$.\\

The trivial but important observation in this setting is that the inclusion map $z : \s \hookrightarrow \hat{\C}$ is automatically equivariant with respect to $\rho = \Id : \Gamma \hookrightarrow \GLtC$, so reversing as it were the independent and dependent variables in the $\D$ operator, Corollary \ref{bijection} can be applied to $f(z) = z$ to generate equivariant rational expressions in automorphic forms and their derivatives. The simplest of these is the first-order expression $\phi_\alpha(z) := \D_{\theta}(z)$ in the coefficient $\alpha \in \M_2$ of an invariant 1-form $\theta = \alpha(z)dz$, $$\phi_\alpha(z) = z + 2 \frac{\alpha(z)}{\alpha'(z)}.$$ We call this the restricted \emph{weight $2$ $\phi$-operator}
\begin{align*}
\phi : \M_2 & \, \to \, \equi\\
\alpha & \; \mapsto \; \phi_\alpha
\end{align*}
which clearly applies unchanged to automorphic forms of non-trivial character:

\begin{example} Given a convex polyhedron $P \subset \R^3$ with barycenter at the origin, consider its radial projection onto the unit sphere and let $v$, $e$, and $f$ be the monic polynomials with zeros at finite stereographic projections of the $V$ vertices, $E$ edge-midpoints, and $F$ face-barycenters, respectively. If each of these sets is an orbit of the double cover $\Gamma$ of the rotation group of $P$, then $v,f,e$ are $\Gamma$-invariants of weights $-V, -F, -E$, respectively, which generate the ring of invariants algebraically, $\mathcal{P}_* \simeq \C[v,f,e]$. Thus $\alpha = \frac{e}{vf} \in \M_2$ and $$\phi_\alpha(z) = z - 2\,\big(\frac{v'}{v} - \frac{e'}{e} + \frac{f'}{f}\big)^{-1}$$ is a degree $|\Gamma|+1$ rational equivariant. Evidentally the result is invariant under duality $v \leftrightarrow f$, up to which each such polyhedron is given by a Platonic solid $P = P_{\{3,n\}}$ with Schl\"afli symbol $\{3,n\}$. Explicit computation of $\phi_\alpha$ is then facilitated by the normalization rotating one vertex onto the positive vertical axis (the point at infinity) and a second onto the positive real axis, so that the invariants $v_n,e_n,f_n$ assume the particularly simple, well-known form appearing in Table \ref{PlatInv}.
\begin{table}[h] \caption{Generators of the invariant ring for normalized Platonic groups}
\centering
\begin{tabular}{|c|l|l|} \hline & & \\[-1.5ex]
$\Gamma$ & $v_n$ vertex invariant,\; $f_n$ face invariant,\; $e_n$ edge invariant & syzygy \\[0.5ex] \hline & & \\[-1.5ex]
$A_4$ & $v_3 = z^3 + \sqrt{2}/4,\;\; f_3 = z^4 - 2\sqrt{2} z,\;\; e_3 = z^6 + 5\sqrt{2} z^3 - 1$ & $e_3^2 - f_3^3 = 16\sqrt{2} v_3^3$ \\[0.5ex] \hline & & \\[-1.5ex]
$S_4$ & $v_4 = z^5 - z,\;\; f_4 = z^8 + 14z^4 + 1,\;\; e_4 = z^{12} - 33z^8 - 33z^4 + 1$ & $e_4^2 - f_4^3 = -108 v_4^4$ \\[0.5ex] \hline & & \\[-1.5ex]
$A_5$ & $v_5 = z^{11} + 11z^6 - z,\;\;\, f_5 = z^{20} - 228z^{15} + 494z^{10} + 228z^5 + 1$ & $e_5^2 - f_5^3 = 1728 v_5^5$\\[0.5ex] & $\quad e_5 = z^{30} + 522z^{25} - 10005z^{20} - 10005z^{10} - 522z^5 + 1 $ & \\[0.5ex] \hline
\end{tabular} \label{PlatInv} \end{table}
The rotation $R \in \SOtR$ restoring $P$ to its original configuration induces M\"obius conjugation $\equi \to T\circ\equi\circ T^{-1}$ on the monoid of equivariants by either pre-image $T \in \SUt$ in the spin cover $\SUt \to \SOtR$.
\end{example}

The resemblance $\phi_\alpha$ bears to Newton's method $N_\alpha = z - \alpha/\alpha'$ is not superficial: The divisor $(\alpha)$ consists of fixed points of both $\phi_\alpha$ and $N_\alpha$, with the important difference that the superattracting case for $N_\alpha$ occurs at simple zeros of $\alpha$, and for $\phi_\alpha$ at double poles of $\alpha$. This coincidence becomes more meaningful for $\phi$ defined by an \emph{exact} 1-form $\theta = df$, giving a second-order operator $f \mapsto \phi_{f'}$ on $\M_0 = \invar \to \equi$, the differential of which satisfies the relation $(\varphi_z f)^2d\phi_{f'} = \S_z f$. Invariance of this relation under inversion $f \mapsto 1/f$ suggests the geometric significance of this special case of $\phi$, which appears implicitly in the next most classical construction in numerical analysis:

\begin{example} Given a function $f \in \mero$, define lifts $$A : \s \to \C^* \ltimes \C, \quad M : \s \to \PGLtC$$ where $A_z \in \aut(\C)$ and $M_z \in \aut(\cpo)$ osculate $f$ at $z$ (that is, the $1$-jet of $A_z$ and $2$-jet of $M_z$ agree with those of $f$ at $z$). Then the rational maps $$N_f(z) = A_z^{-1}(0),\quad H_f(z) = M_z^{-1}(0)$$ define Newton's method and ``projectively natural Newton's method'' \cite{DM89}, better known in numerical analysis as \emph{Halley's method}. The manifest equivariance of $H_f$ with respect to the invariant subgroup of $f$ can also be inferred setting $\theta = d(1/f)$ and factoring $M$ and its group-theoretic inverse $M^{-1}$ using Legendrian curves, $$M = L_{f,dz}L_{z,dz}^{-1}, \quad  M^{-1} = L_{z,\theta}L_{1/f,\theta}^{-1}N,$$ so that $\pi_2 \circ M^{-1} = \pi_2 \circ L_{z,\theta}$, and thus $H_f$ coincides with the $\phi$-operator of $\alpha = f'/f^2$, $$H_f(z) = \D_\theta(z) = z + \frac{2f f'}{f f'' - 2{f'}^2}.$$ In particular, $H_f \in \equi$ if and only if $\alpha \in \M_2$, if and only if $f \in \invar$. A simple zero $z_0$ of $f$ is a double pole of $\alpha$ and thus a second-order superattracting fixed point $H_f(z_0) = z_0$, $H_f'(z_0) = H_f''(z_0) = 0$, which is of third order if and only if $M_{z_0}$ hyperosculates $f$ at $z_0$. %\footnote{Its reputation as ``the most frequently rediscovered method'' in numerical analysis \cite{Tr64} seems apt.}
\end{example}

Forming absolute invariants $h$ as rational expressions of brackets of automorphic forms yields, in the notation of Theorem \ref{operator}, a hierarchy of higher-order extensions $(\Phi\,h)(z)$ of the $\phi$-operator, starting with the following first-order family: Given $\Gamma$-automorphic forms $\alpha$ and $\beta$ of weights $k$ and $k+2$, both $[\alpha,\beta]_1$ and $\beta^2$ are weight $2k+4$, yielding a 2-parameter family of invariants $$h = \lambda_1 + \lambda_2\frac{[\alpha,\beta]_1}{\beta^2},\;\; \theta = (\beta/\alpha)dz$$ independent of the characters of $\alpha$ and $\beta$ (provided they coincide) with corresponding equivariant functions $(\Phi\,h)(z)$ parameterized by constants $\lambda_1, \lambda_2 \in \C$, all of order at most $1$ in both $\alpha$ and $\beta$. The dependence of $(\Phi\,h)(z)$ on $\beta$ and $\beta'$ is eliminated by setting $\lambda_1 = 0$ and $\lambda_2 = -1/2k$, in which case $\phi_{\alpha} := (\Phi\,h)(z)$ defines the \emph{weight $k$ $\phi$-operator} $\M_k \to \equi$, $$\phi_{\alpha} = z + k \frac{\alpha}{\alpha'},$$ which in turn extends to the \emph{bi}weight $(k,k+2)$ $\phi$-\emph{operator} $$\phi : \M_k \times \M_{k+2} \to \equi$$ of order zero in the second factor, $$\phi_{\alpha,\beta} = z + \frac{k\alpha}{\alpha' + \beta},$$ obtained by setting $\lambda_1 = 1/k$ and $\lambda_2/2 = -1/2k$ in $(\Phi\,h)(z)$.\\

That the weight 2 operator on a quotient can be rewritten $\phi_{\beta/\alpha} = \phi_{\alpha\beta, [\alpha,\beta]_1}$ is the first indication that $\equi$ can often be described completely by applying $\phi$ to pairs of \emph{holomorphic} invariants. For example, in the rational category, the restriction $$\phi : \mathcal{P}_{*-1} \times \mathcal{P}_{*+1} \to \equicpo$$ to polynomial invariants is surjective by a construction of Klein \cite{Kl77}:

\begin{example} Given $\Gamma \subset \GLtC$ finite, let $\mathfrak{X}^\Gamma\C^2$ and $\mathfrak{X}_\Gamma\C^2$ denote the linear spaces of $\Gamma$-\emph{invariant} and $\Gamma$-\emph{equivariant polynomial vector fields} on $\C^2$, $$\mathfrak{X}^\Gamma\C^2 = \{X : \C^2 \to \C^2 \,|\, X \circ A = \chi\, X\}, \;\;\; \mathfrak{X}_\Gamma\C^2 = \{X : \C^2 \to \C^2 \,|\, X \circ A = \chi\, A\,X\}$$ with $\chi \in \C[z_1,z_2]$. Then given homogeneous invariant polynomials $A,B \in \C[z_1,z_2]^\Gamma$, Klein \cite{Kl77} defines equivariant polynomial vector fields $X,Y \in \mathfrak{X}_\Gamma\C^2$ by $$X = \nabla_\omega A,\;\; Y = -B\,\E,$$ where $\nabla_\omega = J \circ \nabla = (\partial_{z_2}, - \partial_{z_1})^t$ is the symplectic gradient with respect to $\omega = dz_1\wedge dz_2$ and $\E = (z_1, z_2)^t$ the radial vector field. If $X$ and $Y$ have equal degree, their sum $Z$ is also homogeneous and equivariant. Since the isomorphism $v \mapsto \omega(v,\cdot)$ between $\C^2$ and its dual identifies the symplectic gradient with the exterior differential and equivariant vector fields with \emph{in}variant 1-forms, $\omega(Z,\cdot) = dA - B\lambda$ is invariant and homogeneous, where $\lambda = \half (z_1 dz_2 - z_2 dz_1)$. Doyle and McMullen \cite{DM89} show that any such 1-form can be so expressed,\footnote{Klein arrives at this conclusion as well, citing results of Gordan and Clebsch; see the discussion starting on page 345 of \cite{Kl77}.} thus any equivariant homogeneous vector field is of the form $Z = \nabla_\omega A -B\,\E$. It is then easily verified that if $A, B$ are homogenizations of invariant polynomials $\alpha(z), \beta(z) \in \mathcal{P}_*$, the rational equivariant $[Z] : \hat{\C} \to \hat{\C}$ to which $Z$ descends coincides with their $\phi$-operator, $$\phi_{\alpha,\beta} = [\nabla_\omega A - B\,\E].$$ \end{example}

%\begin{cor}\textnormal{(\cite{Kl77}, \cite{DM89})} Given $\Gamma \subset \GLtC$ finite, the restriction $$\phi : \{(\alpha,\beta) \, | \, \alpha \in \mathcal{P}_k,  \beta \in \mathcal{P}_{k+2}\} \to \equicpo$$ of the $\phi$-operator to polynomial invariants is surjective. \end{cor}

This surjectivity gives information about the structure of $\equicpo$ with respect to the natural stratification of $\C(z) = \bigcup_{d \in \N} \textnormal{Rat}_d$ by degree: $\phi$ maps a bidegree $(d+1,d-1)$ polynomial pair $(\alpha,\beta)$ to a degree $d$ rational function $f = \phi_{\alpha,\beta} \in \textnormal{Rat}_d$, and an arbitrary degree $d$ element $f \in \equicpo$ admits a degree-preserving deformation in $\equicpo$ if and only if there exist polynomials $\alpha$, $\beta \ne 0$ and a constant $\lambda \in \C$ such that $f = \phi_{\alpha,\lambda\beta}$, the dimension of which (as a subvariety of $\textnormal{Rat}_d$) can be determined from the graded dimension of the subalgebra $\mathcal{P}_* \subset \mathcal{M}_*$ by weight (see \cite{MSW17} for a highly detailed, algebraic treatment). As for equivariants themselves, clearly zeros of $\alpha$ are again fixed points of $\phi_{\alpha,\beta}$, while the critical points of $\phi_\alpha$ are determined by the Hessian of the homogenization of $\alpha$,
$$-(k+1)(\phi_\alpha')(\alpha')^2 = [\alpha,\alpha]_2.$$ This is particularly relevant to the dynamics of $\phi_\alpha : \cpo \to \cpo$ when $\alpha$ is a polynomial vanishing along a special orbit of $\Gamma$, in which cases one can show by studying a construction of Doyle and McMullen \cite{DM89} that the critical points form superattracting 2-cycles: %: Averaging a Hermitian form on $\C^2$ to regard $\Gamma \subset \SUt \to \SOtR$ acting on the sphere $S^2 \subset \R^3$

\begin{example} Given a tiling $T = \bigcup_{i=1}^d T_i$ of the unit sphere $S^2 \subset \R^3$ symmetric about the origin $T = -T$ with a transitive rotation group action, let $DM_{T_k} : T_k \to \bigcup_{i\ne k} T_i$ be the conformal map from a tile $T_k$ to the complement of its opposite $-T_k$ fixing its center and mapping each vertex to its antipode. Conjugating with stereographic projection, these glue along edges to form a well-defined, manifestly equivariant rational map $$DM_T : \hat{\C} \to \hat{\C}$$ of degree $d-1$, with fixed points at face-centers $f_T = 0$ and critical points at vertices $v_T = 0$, ramifying to order $\delta - 2$ at a vertex of degree $\delta$. Then $f_T$ and $v_T$ are invariants of weight $-d$ and $-\sum (\delta - 2) = 2V-2E = -2d+4$, respectively. Radial projections of suitable polyhedra in $\R^3$ provide examples of such tilings, including an infinite family of bipyramids or trapezohedra (of dihedral symmetry $\Gamma \simeq D_n$ for $n$ even or odd, respectively), any non-self-dual Platonic solid $P$ (of $\Gamma \simeq S_4$ or $A_5$ symmetry), and any Catalan solid $xP$ obtained by applying the Conway \emph{join}, \emph{ortho}, or \emph{kis} operator to $P$; the rhombic tiling $jP$ replacing the edges of $P$ with those joining vertices to face-centers, the deltoidal tiling $oP$ connecting the edges of $P$ orthogonally to face-centers, or the triangulation $kP$ connecting the vertices of $P$ to face-centers, $$v_{jP} = v^{q-2}f^{p-2},\;\; v_{oP} = v^{q-2}f^{p-2}e^2,\;\; v_{kP} = v^{2q-2}f^{p-2}$$ where $v,f,e$ are the invariants of $P$ and $\{p,q\}$ its Schl\"afli symbol. In each case the bracket $[f_T,f_T]_2$ is proportional to the vertex invariant $v_T$, so by degree, fixed point, and ramification considerations, $DM_T$ coincides with the $\phi$-operator of the face invariant $f_T$, $$DM_T = \phi_{f_T}.$$ All six $A_5$-equivariants so produced are degree-rigid in $\equicpo$, including Klein's example $K = DM_D = \phi_{v_5}$ determined by the dodecahedron\footnote{That the joint basin of attraction of the $20$ critical fixed points $v_D = 0$ of $K \circ K$ has full-measure is an essential feature in the construction of the iterative algorithm appearing in \cite{DM89}.} $D$. Only $DM_C$ and $DM_{jC}$ determined by the cube $C$ are rigid as $S_4$-equivariants, but none is rigid as an $A_4$-equivariant since face diagonals of $C$ form edges of a dual pair of tetrahedra with invariants $v$, $f$, $e$, so that $DM_C = \phi_{e}$ and $DM_{jC} = \phi_{e_4}$ admit $2$-parameter deformations by $\beta = \lambda_1 v + \lambda_2 f \in \mathcal{P}_{-4}$ and $\tilde{\beta} = e\,\beta  \in \mathcal{P}_{-10}$, respectively.
\end{example}

Perhaps unaware of Klein's construction but inspired by that of Heins, Smart \cite{Sm72} rediscovered the weight $k$ $\phi$-operator in the general case of a Kleinian group $\Gamma$ with domain of discontinuity $\s \subset \hat{\C}$, and established the existence of two linearly independent $\Gamma$-automorphic forms $\alpha_0$, $\alpha_1$, $\phi$-transforms $f_0$, $f_1$ of which, together with $f_\infty = z$, determine the inverse $f \mapsto [f_\infty,f_0,f_1,f]$ of a Cayley-type map from $\equi$ to $\invar$. Since Theorem \ref{operator} is an analogue of this map at the level of operators, Corollary \ref{bijection} is its natural generalization on arbitrary surfaces, replacing $z$ with any $f_\infty \in \equi$ and setting $f_0 = \Phi(0) f_\infty$, $f_1 = \Phi(1) f_\infty.$  Brady had come by these in the modular case $\Gamma = \PSLtZ$ by applying the Heins construction to ``pseudo-periodic'' generalizations of the Weierstrass zeta function, where equivariance implies interesting recursions relations among Fourier coefficients:

%Interestingly, as Smart observes, the prototypical modular equivariants appearing in the literature already lie in the image of the weight $k$ operator:

\begin{example} Let $E_{2k}$ denote the weight $2k$ (normalized) Eisenstein series $$E_{2k}(\tau)\; = \frac{1}{2\zeta(2k)}\!\sum_{(m,n) \in \Z^2 - 0}\! (m\tau + n)^{-k}.$$ Writing $\eta_{\Lambda}(\tau)$ and $\eta_{\Lambda}(1)$ in terms of $E_2$, the Heins function $H(\tau) = \frac{\eta_{\Lambda}(\tau)}{\eta_{\Lambda}(1)}$ coincides with that studied by Nahm, $\tau + \frac{6}{\pi i} \frac{1}{E_2},$ equivariance of which follows from $\eta' = \textstyle\frac{i\pi}{12} E_2\eta$ where $$\eta(\tau) = q^{\frac{1}{24}}\prod\limits_{n=1}^{\infty}(1-q^n), \;\; q = e^{2\pi i\tau}$$ is the Dedekind eta function, a modular form of weight $\frac{1}{2}$, so that $H(\tau) = \phi_{\eta}(\tau)$. More generally, since any absolute $\PSLtZ$-invariant $h$ is a rational expression $r \circ j$ of Klein's $j$-invariant, any non-trivial $\PSLtZ$-equivariant $f$ is of the form $\Phi_{X}(r\circ j)\tau$ with $X = d/dj$, expressible in terms of Eisenstein series by applying the Ramanujan identities $$E_2' = 2\pi i\left(\frac{E_2^2 - E_4}{12}\right), \quad E_4' = 2\pi i \left(\frac{E_2E_4 - E_6}{3}\right), \quad E_6' = 2\pi i \left(\frac{E_2E_6 - E_4^2}{2}\right)$$ to differentiate the relation $\dot{\tau}^{-1} = j' = -2\pi i (E_6/E_4)j$, so that
$$f(\tau) \;=\; \Phi_{X}(r\circ j)(\tau) \;=\; \tau + \frac{6}{\pi i}\frac{E_4E_6}{E_2E_4E_6 - 3E_4^3 + (r\circ j)E_6^2} \;=\; \tau + \frac{6}{\pi i}\frac{1}{E_2 + \mu}$$ where $\mu \in \M_2$. This determines a bijection $\M_2 \;\longleftrightarrow\; \mero_{\PSLtZ}(\H^2) - \{\tau\}$, which in turn implies that the restriction $\phi : \mathcal{H}_{*-1} \times \mathcal{H}_{*+1} \to \mero_{\PSLtZ}(\H)$ to (holomorphic)  modular forms $\mathcal{H}_{*} \simeq \C[E_4,E_6]$ is surjective. See Sebbar for a series of detailed studies \cite{SS12}.
\end{example}

When the representation $\rho : \Gamma \to \PGLtC$ is not M\"obius conjugate to the inclusion map $\rho \ne T \circ \iota \circ T^{-1}$, or $\s$ not conformal to a subset of $\cpo$, the classical techniques do not directly assist in constructing an initial equivariant to which Corollary \ref{bijection} can be applied.
If however $\s$ is given as a quotient of such a subset, equivariants can often be obtained by descent. For example, for $\Sigma$ compact of genus greater than $1$ uniformized $\H^2 \to \s \simeq \H^2/\Gamma'$ by a Fuchsian group $\Gamma' \subset \PSLtR$, any non-trivial projective representation of a larger group $\rho : \Gamma \to \PGLtC$ with $\Gamma' \subset \ker\rho$ descends to the quotient $G = \Gamma/\Gamma'$, with respect to which a $\Gamma$-equivariant function on $\H^2$ descends to an element of $\equis$. In genus zero, this logic can be reversed: Given a Fuchsian group $\Gamma \subset \SLtR$ admitting a finite quotient $G = \Gamma/\Gamma'$ by a normal subgroup $\Gamma'$ of genus zero (that is, the compactification $\s = \overline{\H^2/\Gamma'}$ has genus zero), the field of $\Gamma'$-invariants on $\H^2$ is generated by a single function, a \emph{Hauptmodul} $f \in \mero^{\Gamma'}$, unique up to M\"obius transformation, such that $\mero^{\Gamma'} = \C(f)$. Since any composition $f \circ g$ with $g \in \Gamma$ is again $\Gamma'$-invariant, there exists a rational function $\rho(g) \in \C(z)$ such that $f \circ g = \rho(g) \circ f$, and since suitable compositions of these $\rho(g)$ yield the identity, they must be linear fractional, defining a representation $\rho : \Gamma \to \PGLtC$ with respect to which $f$ is tautologically equivariant:
%The rational case becomes relevant in this context since $G \subset N(\Gamma')/\Gamma' \subset \auts$ is finite, so that $\equi$ carries a monoid action of $\mero_G(\cpo)$ by post-composition.
%...of such quotients, and is implicit in certain number theoretic applications
%, the Hauptmoduln for which are classically known in terms of eta quotients and infinite $q$-products

\begin{example}
The \emph{principal congruence subgroup} given by the kernel $\Gamma(n) = \ker \pi_n$ of the projection $\pi_n : \, \SLtZ \to \SLtZn$ from $\Gamma = \SLtZ$ is genus zero for $1\le n \le 5$, and the representation $\rho_n : \Gamma \to \PGLtC$ induced by a choice of Hauptmodul $j_n$ can be explicitly determined by the generators $$S : \tau \mapsto -1/\tau, \quad T : \tau \mapsto \tau + 1$$ of $\Gamma$. The kernel of the latter is again $\Gamma(n)$, hence the image $\rho_n(\Gamma) \simeq \Gamma/\Gamma(n)$ is trivial for $n=1$, where $j_1 = j$ is the Klein $j$-invariant, the anharmonic group $S_3$ for $n=2$, where $j_2 = \lambda : \H \to \C - \{0,1\}$ is the modular lambda function, and a group of Platonic symmetries otherwise, and therefore a full solution of Problem \ref{p1} with respect to $\rho_n$ is $$\mero_{\Gamma}(\H^2) = \mero_{\rho_n(\Gamma)}(\cpo) \circ j_n.$$ More generally, there are $132$ distinct conjugacy classes of genus zero congruence subgroups (that is, subgroups $\Gamma' \subset \SLtZ$ containing some $\Gamma(n) \subset \Gamma'$), including $$\Gamma_1(n) = \{A \in \SLtZ \;|\; \pi_n(A) = \begin{pmatrix} 1 & * \\ 0 & 1 \end{pmatrix} \}$$ for $n=1,2,\ldots,10$ and $12$, and the \emph{Hecke congruence subgroups} $$\Gamma_0(n) = \{A \in \SLtZ \;|\; \pi_n(A) = \begin{pmatrix} * & * \\ 0 & * \end{pmatrix} \}$$ for $n=1,2,\ldots,10$, $12,13,16,18,$ and $25$. The former is a normal subgroup of the latter with quotient $\Gamma_0(n)/\Gamma_1(n) \simeq (\Z/n\Z)^\times$, while the full normalizer of each in $\PSLtR$ can be found in Lang \cite{La01} and Akbas and Singer \cite{AS90}, respectively. For any such group $\Gamma'$, since $j$ is also $\Gamma'$-invariant it is rationally related $j = r \circ f$ to any $\Gamma'$-Hauptmodul $f$,
\begin{center}\begin{tikzcd}[column sep=tiny]
\H \arrow{dr}[swap]{j}\arrow{rr}{f} & & \C \arrow{dl}[dashed]{r} \\
& \C &
\end{tikzcd}\end{center}
by an absolute invariant $r \in \mero^{\rho(\Gamma)}(\cpo)$. This invariant is given explicitly by $$r_n(z) = k_n f_n(z)^3/v_n(z)^n$$ for the $j_n$ of Table \ref{Haupt}, $j = r_n \circ j_n$, with $v_n$, $f_n$ listed in Table \ref{PlatInv} together with the $S_3$-invariants $f_2(z) = z^2 - z +1$ and $v_2(z) = z^2-z$, and suitable constants $k_n$.
\begin{table}[h] \caption{Hauptmoduln for principal congruence subgroups}%, $[N] := \eta(N\tau)$}
\centering
\begin{tabular}{|c|c|c|l|} \hline & & & \\[-1.5ex]
$\Gamma'$ & $\Gamma/\Gamma'$ & $k_n$ & Hauptmodul $j_n(\tau)$, $q = e^{2\pi i\tau}$ \\[0.5ex] \hline & & & \\[-1.5ex]
$\Gamma(2)$ & $S_3$ & $256$ & $j_2(\tau) = \eta(\tau/2)^8\eta(2\tau)^{16}/\eta(\tau)^{24}$  \\[0.5ex] \hline & & & \\[-1.5ex]
$\Gamma(3)$ & $A_4$ & $-54\sqrt{2}$ & $j_3(\tau) = -\sqrt{2}\eta^3(\tau/3)/6\eta^3(3\tau)-2\sqrt{2}$ \\[0.5ex] \hline & & & \\[-1.5ex] % = [1/3]^3/[3]^3
$\Gamma(4)$ & $S_4$ & $16$ & $j_4(\tau) = 2\eta^2(\tau)\eta^4(4\tau)/\eta^6(2\tau) = 2q^{1/4}\prod_{n\ge 1}(1-q^{2n-1})/(1-q^{4n-2})^2$ \\[0.5ex] \hline & & & \\[-1.5ex] % = 2q^{1/4}\prod_{n\ge 1}(1+q^{n})^{(-1)^n 2} = [1]^2[4]^4/[2]^6
$\Gamma(5)$ & $A_5$ & $-1$ & $j_5(\tau) = q^{1/5}\prod_{n\ge 1}(1-q^{n})^{(\frac{5}{n})}$,\; $(\frac{\cdot}{\cdot})$ the Kronecker symbol \\[0.5ex] \hline
\end{tabular}
\label{Haupt}
\end{table}
The equivariance of congruence Hauptmoduln, combined with their relationship to the $j$-invariant and basic results from the theory of complex multiplication, has interesting number theoretic consequences, such as those described by Duke \cite{Du05} concerning values of several related continued fractions due to Ramanujan, particularly the \emph{Rogers-Ramanujan continued fraction} $$R(\tau) = \cfrac{q^{1/5}}{1+\cfrac{q}{1+\cfrac{q^2}{1+\cfrac{q^3}{1+\ddots}}}}$$ which, remarkably, coincides with $j_5(\tau)$.
\end{example}

\begin{example} Consider the homomorphism $\rho : \Gamma(2) = \langle ST^2S,T^2 \rangle \to \Gamma(1)$ determined by $$\rho(ST^2S) = ST^{-1}, \quad \rho(T^2) = S,$$ the kernel of which is the normalizer of $\langle ST^6S,T^4 \rangle$ in $\Gamma(2)$. Kaneko and Yoshida \cite{KY03} construct the \emph{kappa function} $\kappa(\tau)$ as the $\rho$-equivariant map such that the modular lambda $\lambda = J \circ \kappa$ factors through the normalized Klein invariant $J = j/1728$,
\begin{center}
\begin{tikzcd}[column sep=tiny]
\H \arrow{dr}[swap]{\lambda}\arrow{rr}[dashed]{\kappa} & & \H \arrow{dl}{J} \\%- \Gamma(1){i}
& \C & % - \{0,1\}
\end{tikzcd}
\end{center} %$$\mero_{\Gamma(2)}(\H) \;=\; \{\kappa\} \,\cup\, \{\Phi_{X}^{r(\lambda)}\kappa\; | \; r(z) \in \C(z) \}.$$ $$\Phi_{d/d\lambda}^{r(\lambda)}\kappa = \kappa + \frac{J'(\kappa)}{r(\lambda)J'(\kappa)^2 + \half J''(\kappa)} = (\tau + \frac{\dot{\tau}}{r(J)- \half \dot{\tau}^{-1}\ddot{\tau}}) \circ \kappa = (\Phi_{d/dJ}^{r(J)}\tau)\circ\kappa,$$
so by Corollary \ref{bijection}, $\mero_{\Gamma(2)}(\H) - \{\kappa\} = \{\Phi_{X}(r\circ\lambda)\kappa\; | \; r \in \C(z) \}$ where $X = d/d\lambda$. The change of variables formula for the $\D$ operator then implies
$\Phi_{X}(r\circ\lambda)\kappa = (\Phi_{Y}(r\circ J)\tau) \circ \kappa$ with $Y = d/dJ$,
thus $\mero_{\Gamma(2)}(\H) = \mero_{\Gamma(1)}(\H^2)\circ \kappa$. See Tanaka \cite{Ta07} for several related triangle groups $\tilde{\Gamma}$ with morphisms $\tilde{\rho} : \tilde{\Gamma} \to \Gamma(1)$ and explicit corresponding $\tilde{\rho}$-equivariant versions of the kappa function.
\end{example}

\section{Non-commutative function algebras in several variables}\label{SectionNC}

The real virtue of Theorem \ref{operator} is the ease with which it extends to a more general context: Replacing the curve $\s$ with a complex manifold $M$ and the plane $\C$ with a unital complex Banach algebra $\A$, let $\meroA = \mero \otimes_\C \A$ be the sheaf algebra of $\A$-valued meromorphic functions on $M$. Denoting the multiplicative identity element $\1 \in \A$ and Banach-Lie group of units $\units \! \subset \A$, there is a central copy of the function field $\meroM \cdot \1 \subset \meroA$ and a pseudo-action on $\A$ by the group $\GLtA = \aut(\A^2)$ of (right $\A$-module) automorphisms of $\A^2 = \A \times \A$, or its projectivization $\PGLtA = \GLtA/\mathcal{Z}(\units)$, by \emph{generalized M\"obius transformations}, defined for $\z \in \A$ and (representative) $T \in \GLtA$ by $$\z \, \mapsto \, T \circ \z = \begin{pmatrix} \a & \b \\ \c & \d \end{pmatrix} \circ \z \, = \, (\a\z+\b)(\c\z+\d)^{-1}$$ if $\c\z+\d \in \units$. The natural generalization of Problem \ref{p1} is then: Given a subgroup $\Gamma \subset \autM$ and projective representation $\rho : \Gamma \to \PGLtA$, describe the set $\equiA$ of $f \in \meroA$ such that $f \circ A = \rho(A) \circ f$ for all $A \in \Gamma$.  \\

To formulate a corresponding problem for differential operators, we must first extend their definition to maintain an unambiguous pseudo-action on $\meroA$. We distinguish two cases: If $\A$ is commutative, a \emph{commutative rational operator} means as before an element of the multiplicative fraction field $\op$ of the ring of meromorphic differential operators on $M$, while if $\A$ is not commutative, a larger set $\opnc$ of non-commutative operators is required. For this we borrow from the treatment of non-commutative function algebras appearing in Kaliuzhnyi-Verbovetskyi and Vinnikov \cite{KV14}: Given a sequence $X = \{X_k\}_{k=0}^\infty$ of differential operators on $M$, we associate to a list of natural numbers $\kappa = (k_0,k_1,\ldots\ k_m)$ a non-commutative monomial $X_\kappa = X_{k_0}X_{k_1}\cdots X_{k_m}$, so that given a subring $R \subset \mero$, a \emph{non-commutative polynomial} $P$ \emph{in $X$ over} $R$ is a finite sum $$P(X) = P(X_0,X_1,\ldots,X_n) = \sum a_\kappa X_\kappa$$ where $a_\kappa \in R$ and $n = \max\{n\in\kappa \,|\, a_\kappa\ne0\}$. Endowed with the obvious notions of addition and multiplication, the set of such polynomials is then isomorphic to the free associative algebra generated by $X$ over $R$. A \emph{non-commutative rational expression} $E$ is now any obtained by successive application of addition, multiplication, and formal multiplicative inverse $P \mapsto P^{-1}$ operations to a finite set of polynomials. Its \emph{domain} $\textnormal{dom}\,E \subset \meroA$ consists of those elements $f \in \meroA$ such that all inversions appearing in the evaluation $E(Xf) = E(X_0f,X_1f,\ldots,X_nf)$ exist on the complement of a subset in $M$ of codimension at least $1$, and two expressions $E_1$, $E_2$ are \emph{equivalent} if they agree on their common (non-empty) domain: $$E_1 \sim E_1 \;\; \Longleftrightarrow \;\; E_1(Xf) = E_2(Xf)\;\; \forall\, f \in \textnormal{dom}\,E_1 \cap \textnormal{dom}\,E_2 \ne \emptyset.$$ A \emph{non-commutative rational operator $\r$ in $X$ over $R$} is then an equivalence class $\r = [E]$ with $\textnormal{dom}\,E \ne \emptyset$, the set $R\X$, or simply $\X$, of which is again an algebra, isomorphic to the free skew-field generated by $X$, as well as a monoid with respect to operator composition. In the sequel, a meromorphic vector field $X$ will also denote its set of iterates $X_k = X \circ X \circ \cdots \circ X$, with $X_0 = \Id$, $R = \mero$ unless explicitly stated otherwise, and $\opnc$ the algebra of sums, products, and compositions of operators generated from all vector fields on $M$. Since such an operator has no preferred representative expression, there is no natural inclusion $\op \hookrightarrow \opnc$, but a section $\mathfrak{U} \subset \op \to \opnc$ of natural projection to the quotient by the multiplicative commutator (the \emph{classical limit}) $$\pi : \opnc \to \opnc/[\opnc,\opnc] \simeq \op$$ will be called a \emph{quantization} of $\mathfrak{U}$, and the task of selecting a quantization with prescribed properties the \emph{ordering ambiguity problem} on $\mathfrak{U}$. \\
%Clearly $\X \subset \opnc$ is a submonoid and depends only on the span of $X$ in $TM$.
%\footnote{We also include the constant operator $\1 : f \to \1$ among the generators by convention, evidentally equivalent to expressions of the form $E_1 = EE^{-1}$ or $E_2 = E^{-1}E$ for any $E$ with $\textnormal{dom}\,E \ne \emptyset$.}
%Note that such an operator $\r \in \opnc$ may admit many formally different but algebraically equivalent expressions, and not necessarily any unique, preferred simplification, but still has an \emph{order} (the minimum over all its algebraic simplifications of the maximal number of compositions of vector fields occurring in any term) as well as a \emph{depth} (the minimum over all simplifications of the maximal number of nested inverse operations appearing in any term).

With this larger monoid in place, we seek to describe the submonoid $\opnceq \subset \opnc$ of \emph{equivariant} rational operators, consisting of $\r \in \opnc$ such that $\r \circ T = T \circ \r$ for all $T \in \PGLtA.$  As we will see, this amounts to an ordering ambiguity problem over $\opeq$, and we again look to the geometry of the M\"obius group for its resolution. Observe first that, as in the 1-dimensional case, the pseudo-action of $\PGLtA$ on $\A$ comes from a genuine action on the projective line $\apo$ appearing in quantum mechanics \cite{BN05}: We call a right $\A$-submodule $\ell \subset \A^2$ a \emph{line} if it is isomorphic to $\A$ and admits a complementary right $\A$-submodule $\ell'$ with $\ell \oplus \ell' \simeq \A^2$, and define $\apo$ as the space of lines in $\A^2$. Such a line $\ell = \{ (\x\a,\y\a) \in \A^2 \; | \; \a \in \A \} \in \apo$ can be represented by homogeneous coordinates $[\x,\y]^t$ with $\x, \y \in \A$, or by the affine coordinate $\z = \x\y^{-1} \in \A$ if $\y \in \units$, hence the familiar embedding $\A \hookrightarrow \apo$ given by $\z \mapsto [\z,\1]^t$ and common interpretation of $\apo$ as a ``projective completion'' of $\A$ (see Bertram and Neeb \cite{BN05} and references therein). This embedding permits the identification of a function $f \in \meroA$ with a morphism $[f,\1]^t : M \to \apo$ to which it extends, and exhibits M\"obius transformation on $\A$ as projective action of $\PGLtA$ on $\apo$. Both actions are transitive, as is that of $\PGLtA$ on space of complementary line pairs, since the column vectors of an element $T \in \GLtA$ provide homogeneous coordinates of such a pair, $\pi_1(T) = T[\1,\0]^t$ and $\pi_2(T) = T[\0,\1]^t,$ and conversely. \\ %and we again have projection maps $\pi_i : \PGLtA \to \apo$ defined by  %, so we have a double fibration...\\

An intuitive analogue of the contact condition for maps into $\PGLtA$ is now evident: Decomposing the left-invariant Maurer-Cartan form $\omega$ on $\PGLtA$ into $\A$-valued component 1-forms $$\omega = \begin{pmatrix} \Xi & \hat{\Theta} \\ \Theta & \hat{\Xi} \end{pmatrix} : T\PGLtA \to \gla,$$ we say that a meromorphic map $L : M \to \PGLtA$ is \emph{reduced with respect to} $X \in T_pM$ if the diagonal entries of $L^*\omega(X)$ at $p$ are both the zero element of $\A$ and the lower-left entry is the multiplicative identity, $$L^*\Xi(X) = L^*\hat{\Xi}(X) = \0, \; L^*\Theta(X) = \1.$$ Although uniqueness fails in higher dimension, the next lemma establishes the local existence of reduced lifts of suitably regular functions: We say that a holomorphic map $f : M \to \apo$ is \emph{regular with respect to} $X \in T_pM$ at $p \in M$ if the push-forward $f_*X \in T_{f(p)}\apo \simeq \textnormal{Hom}(\ell,\A^2/\ell)$, viewed as a morphism from the image $\ell = f(p)$ to the quotient $\A^2/\ell$, is invertible. If we let $f$ also denote the affine coordinate of this map, the condition simply says that the directional derivative at $p$, viewed as an element of $\A$, is invertible, $\dot{f}(p) = X_p(f) \in \units$.

\begin{lem}
Given $f \in \meroA$ regular at $p$ with respect to $X \in \merofields$, there exists a lift $L_{f,X} : U \subset M \to \PGLtA$ reduced with respect to $X$ on a neighborhood $U$ of $p$, unique up to a right $\A$-multiple constant along the flow of $X$.

\begin{proof} By the regularity condition, there is a neighborhood $\tilde{U}$ of $p$ on which the tangent vector $\dot{f} : \tilde{U} \to \units$ is everywhere invertible. Then on a possibly smaller neighborhood $U \subset \tilde{U}$ there exists a holomorphic solution $\Psi : U \to \units$ of the first-order gauge equation $$\Psi^*\omega_R(X) = -\half \dot{f}^*\omega_L(X),$$ where $\omega_R$ and $\omega_L$ denote the right- and left-invariant Maurer-Cartan forms on $\units$. Then $L_{11} = f\Psi$ and $L_{21} = \Psi$ give homogeneous coordinates for $f = [L_{11},L_{21}]^t$ which together with $L_{12} = \dot{L}_{11}$ and $L_{22} = \dot{L}_{21}$ constitute entries of a lift $L_{f,X} : U \subset M \to \GLtA$ of $f$ over $U$, easily verified to be point-wise invertible and reduced.\footnote{As our terminology is meant to suggest, this lemma is again merely the end result of the reduction procedure of Cartan's method of moving frames.}  Since $\Psi$ is unique up to a right multiple (constant with respect to $X$), so is $L$.
\end{proof}
\end{lem}

The regular set $M_0$ of a function $f \in \meroA$ with respect to a vector field $X$ can thus be covered by open sets $U_\alpha$ on which there exist reduced lifts $L_{f,X}^{\,\alpha} : U_\alpha \to \PGLtA$
\begin{center}\begin{tikzcd}[row sep=large]
& \PGLtA \arrow{dl}[swap]{\pi_1}\arrow{dr}{\pi_2} & \\
\apo & U_\alpha \subset M_0 \arrow{l}[swap]{f}\arrow[u, dashed, "L_{f,X}^{\,\alpha}" description]\arrow[r, dotted] & \apo
\end{tikzcd}\end{center}

\noindent from which we extract a pair of (\emph{a priori} locally defined) rational operators $$\D_X f = \pi_2 \circ L_{f,X}^{\,\alpha}, \quad \quad \S_X f = (L_{f,X}^{\,\alpha})^*\hat{\Theta}(X).$$ These exhibit the same variance properties their commutative counterparts:

\begin{prop}\label{operators2} $\D$ is projectively equivariant and $\S$ is projectively invariant.
\begin{proof}
It is sufficient to verify the equivariance property $L_{T \circ f,X} = TL_{f,X}$ on the generators of $\PGLtA$ determined by affine transformations of the form $\z \mapsto \a\z\d^{-1} + \b\d^{-1}$ and the involution $\z \mapsto \z^{-1}$. Indeed, if $\Psi$ solves $\Psi^*\omega_R = -\half \dot{f}^*\omega_L$, then we have
\begin{align*}
(\d\Psi)^*\omega_R &=  -\half \Ad_\d(\dot{f}^*\omega_L) = -\half (X(\a f\d^{-1} + \b\d^{-1}))^*\omega_L,\\
(f\Psi)^*\omega_R &= \half (f^{-1}\dot{f}f^{-1})^*\omega_L = -\half (X(f^{-1}))^*\omega_L.
\end{align*}
\end{proof}
\end{prop}

Instances of both operators occur at least implicitly in differential geometry and analysis, such as a matrix incarnation of the $\D$ operator on $\mero(\mathbb{D},\glnC)$ over the disk $\mathbb{D} \subset \C$, which can be viewed as a complexification of an elementary local construction in Riemannian geometry \cite{Fo86} whereby the Levi-Civita connection is recovered from geodesic flow:

\begin{example} Given a real $n$-dimensional Riemannian manifold $(N,g)$ with geodesic flow $\{\phi_t\} \subset \textnormal{Diff}(TN)$ in the tangent bundle $\pi : TN \to N$, let $V_{\phi_t v}$ be the fibers of the vertical subbundle $V = \ker \pi_* \subset TTN$ along the orbit $\phi_t v$ of a point $v \in TN$, and define subspaces $\gamma_v(t) \subset T_vTN$ by $$\gamma_v(t) = (\phi_{-t})_*V_{\phi_t v}.$$ Setting $\A = \glnR$ and choosing a basis of $T_vTN$ to identify $\apo$ with the Grassmannian $\grassnR$ of $n$-planes in $T_vTN \simeq \R^{2n}$, if $\gamma_v = [f,I]^t : M = (-\epsilon,\epsilon) \subset \R \to \apo$ has reduced lift $L(t) = L_{f, X}$ with respect to $X = d/dt$, then the projections of $L(0)$ give the Riemannian splitting of $T_vTN$, so the affine coordinate of the horizontal subspace $H_v \subset T_vTN$ determined by the Levi-Civita connection at $v$ coincides with $\D_X f(0)$, $$T_vTN = \begin{bmatrix} f(0) \\ I \end{bmatrix} \oplus \begin{bmatrix} \D_X f(0) \\ I \end{bmatrix} =  V_v \oplus H_v$$ as soon as these $n$-planes lie in the big cell $\A \subset \apo$.  See \cite{APD09} for a more general approach to this construction.\end{example}

A new feature appearing in the non-commutative setting is that the local nature of $\S$ is imposed by the topology\footnote{This is somewhat unsurprising considering that the classical Schwarzian derivative, properly regarded, takes values not in the function algebra but in the space of sections of a generally non-trivial line bundle.} of the underlying manifold $M$, since on a neighborhood $U \subset M_0$ in the regular set of $f$ with respect to $X$, $\S_X f = \Psi^{-1}\ddot{\Psi}$ is only well-defined up to conjugation by a constant of integration in $\A$, two solutions $\Psi_\alpha$ and $\Psi_\beta$ of the gauge equation differing by a constant right multiple on an overlap $U_\alpha \cap U_\beta$. We therefore consider the following \emph{global} operators of weaker invariance:

\begin{defn} A rational operator $H \in \opnc$ is called \emph{projectively} \emph{semi-invariant} if $$H(T \circ f) = (\c f+\d)\,H(f)\,(\c f+\d)^{-1}$$ for any representative $T\in\GLtA$.
\end{defn}

\noindent The set of these $\PGLtA$-semi-invariants forms an subalgebra $\opncinv \subset \opnc$ which collapses to that of fully invariant operators in the classical limit, $\opncinv \to \opinv$. To construct such operators, factor the local lift $L$ as a product $L = L^0 A$ of a global lift $L^0 : M_0 \to \PGLtA$ and local gauge $A = \textnormal{diag}[\Psi,\Psi] : U \subset M_0 \to \PGLtA$, so that $\S$ decomposes as $$\S_X f = \Psi^{-1}[X \circ \varphi_X + \varphi_X^2]\Psi,$$ where $\varphi_X = -\half X^{-1}X_2$ is the \emph{pre-Schwarzian}. The latter generates rational operators in $\X$ invariant under $f \mapsto \a f + \b$, the \emph{left-affine group} $\units \ltimes \A$, the algebra of which is isomorphic to $\langle\varphi_X\rangle = \langle\varphi_0,\varphi_1,\varphi_2,\ldots\rangle \subset \X$ where $\varphi_0 = \1$ is a constant operator and $$\varphi_k = -\half X^{-1}X_{k+1}$$ for $k\ge 1$. A polynomial subalgebra of $\Z\langle\varphi_X\rangle$ consisting of projective semi-invariants can now be obtained by means of a modified composition operation on $\opnc$ which reduces to standard operator composition in the classical limit:

\begin{defn} The binary operation $\circ_q : \opnc \times \opnc \to \opnc$ defined on generators by $$\Id \circ_q H = H,\;\;\; X \circ_q H = X \circ H - [\varphi_X, H],$$ will be called \emph{q-composition}, where $[\cdot,\cdot]$ denotes the multiplicative commutator bracket.
\end{defn}

\noindent To an ordered partition $\kappa = (k_0, k_1, \ldots, k_m)$ of a positive integer $|\kappa| = k_0 + k_1 + \ldots + k_m$ we again associate the monomial $\varphi_\kappa = \varphi_{k_0}\varphi_{k_1}\cdots\varphi_{k_m}$ and define $\varphi$-polynomials $$S_n = \sum_{|\kappa|=n+1} a_\kappa \varphi_\kappa\; \in\; \Z\langle\varphi_X\rangle$$ by the sequence\footnote{When $\A = \glnC$, $S_1$ coincides with the ``matrix Schwarzian derivative'' appearing in Schwarz \cite{Sc79}, Zelikin \cite{Ze92}, and numerous other contexts despite not being fully projectively invariant (see however \cite{Ov93}).}
\begin{align*}
S_0 &= \1,\\
S_1 &= \varphi_2 + 3\varphi_1^2,\\
S_2 &= \varphi_3 + 4\varphi_2\varphi_1 + 4\varphi_1\varphi_2 + 12\varphi_1^3,\\
S_3 &= \varphi_4 + 5\varphi_3\varphi_1 + 5\varphi_1\varphi_3 + 6\varphi_2^2 + 24\varphi_2\varphi_1^2 + 20 \varphi_1\varphi_2\varphi_1 + 24\varphi_1^2\varphi_2 + 72\varphi_1^4,\;\;\ldots
\end{align*}
obtained inductively according to $S_{n+1} := X \circ_q S_n$ for $n\ge 1$.

\begin{prop}\label{semi} The operators $S_n$ generate the subalgebra $\Xinv := \opncinv \cap \X$ of projective semi-invariants in $X$, that is, $\Xinv = \langle S_X \rangle = \langle S_0, S_1, S_2,\ldots\rangle$.

\begin{proof} Observe that $S_1 = \Psi^{-1}\S \Psi$, and by the proof of Proposition \ref{operators2}, $\Psi$ transforms analogously to a weight $1$ automorphic form, $\Psi(T \circ f) = (\c f + \d)\Psi(f)$ which, together with the invariance of $\S$, implies $S_1 \in \Xinv$. The gauge equation $X \circ \Psi = \varphi_X\Psi$ together with the Leibniz rule implies that if $\mathcal{H} = \Psi^{-1}H\Psi$ is a local operator with $H \in \X$, then $$X \circ \mathcal{H} = \Psi^{-1}(X \circ_q H)\Psi,$$ so the remaining $S_n$ are related to the (fully $\PGLtA$-invariant) higher derivatives $X_n \circ \S$ by $S_{n+1} = \Psi^{-1}(X_n \circ \S_X)\Psi$, and restricting $q$-composition to $\X$ preserves the subalgebra of semi-invariants, $\X \circ_q\Xinv = \Xinv$. By Cartan, the subalgebra $\X \circ \S$ gives the full set of local invariants under the prolonged action of $\PGLtA$ on the sub-bundle generated by $X$ of jets of maps from $M$ to $\A$ which, since conjugation by $\Psi$ defines an isomorphism between this subalgebra and $\Xinv$, completes the proof.
\end{proof}
\end{prop}

Thus on a curve $\s$, both $\opnc \circ \S$ and $\opnc \circ_q S$ provide solutions to the ordering ambiguity problem over $\mathfrak{U} = \opinv$, the former fully invariant but generally only locally defined for $\pi_1(\s)$ non-trivial, the latter global but only semi-invariant for $\A$ non-commutative. The second-order operator $\D_X = \Id + X \varphi_X^{-1}$ however gives a fully equivariant, global quantization of its commutative counterpart without integrating the gauge equation, $$\D_X f = f - 2 \dot{f} \ddot{f}^{-1} \dot{f}.$$ This asymmetry in the nature of $\D$ and $\S$ on non-commutative function spaces is in fact a general phenomenon distinguishing equivariant and invariant operators: The deformation of $\D$ by $H \in \opnc$ via $\Phi_X : \opnc \to \opnc$ formally identical to that of the commutative case $$\Phi_X H := \Id + X[H + \varphi_X]^{-1}$$ is again global, and $\PGLtA$-equivariant if and only if $H$ is $\PGLtA$-\emph{semi}-invariant, a fact which, together with Proposition \ref{semi}, establishes the advertised solution of the ordering problem over $\opeq$ on a curve:

\begin{prop}
Given $X \in \merofields(M)$, $\Phi_X$ restricts to $\opncinv \to \opnceq$ injectively, and further to $\Xinv \to \Xeq - \{\Id\}$ bijectively.
\begin{proof} The transformation properties of $L = L^{0} \textnormal{diag}[\Psi,\Psi]$ and $\Psi$ imply that of $L^{0}$, namely $$L_{T \circ f,X}^{\,0} = T\,L_{f,X}^{\,0}\,\textnormal{diag}[(\c f+\d)^{-1},(\c f+\d)^{-1}],$$ so writing $\Phi_X H(f) = L_{f,X}^{\,0} \circ H(f)$ with $L^0$ acting point-wise implies $\Phi_X(\opncinv) \subset \opnceq$, as well as the already obvious invertibility of $\Phi_X : \opnc \to \Phi_X(\opnc)$ upon restricting to its image, which omits only $\Id$ in $\Xeq$.
\end{proof}
\end{prop}

Theorem \ref{deform} is an immediate corollary: Given $\Gamma \subset \autM$, $\rho : \Gamma \to \PGLtA$, and $X \in \merofieldsinvar$ a $\Gamma$-invariant vector field, an element $f \in \equiA$ with meromorphically invertible derivative $\dot{f} = X(f)$ belongs to the domain $\textnormal{dom}\,H$ of any polynomial $H \in \invar\langle S\rangle$ in $S_n$ over $\invar$, so that $$f+\dot{f}[H(f)-\half\dot{f}^{-1}\ddot{f}]^{-1} \in \equiA$$ is well-defined as soon as $H(f)-\half\dot{f}^{-1}\ddot{f}$ is meromorphically invertible as well. Since $\1+t\dot{f}^{-1}\ddot{f}$ is invertible for $t \in \mathcal{U} \in \invar$ in some neighborhood of the zero function $\0 \in \mathcal{U}$, constant operators define a deformation $\{f_t\}_{t \in \mathcal{U}}$ of $f = f_0$, namely $$f_t = f+t\dot{f}[\1-\textstyle{\frac{t}{2}}\dot{f}^{-1}\ddot{f}]^{-1}.$$ If additionally $\ddot{f}^{-1} \in \meroA$, then $\hat{f} = \D_X f$ admits a larger deformation $\{\hat{f}_H\}_{H \in \mathfrak{U}}$ in a neighborhood $\mathfrak{U} \subset \{H \in \Xinv \,|\, f \in \textnormal{dom}\,H\}$ of the zero operator, $$\hat{f}_H = f+\dot{f}[\dot{f}H(f)-\half\ddot{f}]^{-1}\dot{f}.$$ %with $\hat{f}_H \to f$ as $H \to \infty$.\\

As in the classical setting for curves, many triples $(M,\Gamma,\A)$ admit a ``trivial'' element $Z \in \equiA$, a holomorphic embedding $Z : M \hookrightarrow\apo$ equivariant with respect to an inclusion $\rho : \aut(M) \hookrightarrow \PGLtA$. Such manifolds include all irreducible Hermitian symmetric spaces of types II and III, as well as certain submanifolds of these and their products, realizable as bounded domains in some $\A = \glnC$. Our main interest is the case of $\Gamma$ countably infinite, in which case there typically exist $m = \dim_\C M$ algebraically independent invariant functions. If moreover $\Gamma$ admits a \emph{complete system of factors of automorphy} $\{\eta_{i,T}\}$, that is, functions $\eta_{i,T} \in \mero$ indexed by $i\in\N,T\in \Gamma$ satisfying $\eta_{i,ST}(Z) = \eta_{i,S}(T \circ Z)\,\eta_{i,T}(Z)$, holomorphic $h_i \ne 0$ such that $h_i(Z)/\eta_{i,T}(Z)$ is holomorphic, and subgroups $\Gamma_i \subset \Gamma$ such that $\eta_{i,T} \equiv 1$ for $T \in \Gamma_i$ and\footnote{We have omitted an additional technical condition; see \cite{BG55} for the precise definition.} the sum $$\sum_{T \in \Lambda_i}h_i(Z)/\eta_{i,T}(Z) \ne 0$$ over any set of representatives $\Lambda_i \subset \Gamma$ of right cosets of $\Gamma_i$ converges uniformly on compact subsets, Bochner and Gunning \cite{BG55} establish the existence of $\dim_\C M$ \emph{functionally} independent $z_k = f_{i_k j_k}$ selected from the invariants $f_{ij} = g_{ij}(g_{i1})^{-j}$ where $g_{ij} = \sum_{T\in\Lambda_i}(\eta_{i,T})^{-j}$. The presence of such a system is again typical of classical examples, in which case these $\{z_k\}_{k=1}^m$ generate $\invar$ algebraically and furnish a $\C$-basis $\{\frac{\partial}{\partial z_k}\}_{k=1}^m \subset \merofieldsinvar$ of the holomorphic tangent space away from the divisor of $dz_1 \wedge dz_2 \wedge \ldots \wedge dz_m$. %The simplest non-compact case is again the modular setting:

\begin{example} Let $M = \siegeln$ be the \emph{degree} $n$ \emph{Siegel upper-half space} of symmetric complex matrices with positive definite imaginary part, $\siegeln = \{Z \in \glnC \, | \, Z^t = Z, \, \textnormal{Im}(Z)>0 \}$, and $\Gamma = \SptnZ$ the Siegel modular group. Then given any set of representatives $\Lambda \subset \Gamma$ of right cosets of the subgroup $$\Gamma' = \{\begin{pmatrix} U^t & SU^{-1} \\ 0 & U^{-1} \end{pmatrix} \in \Gamma  \;|\; U \in \SLnZ,\, S^t = S \},$$ the weight $k$ Siegel Eisenstein series $$E_k(Z) = \sum_{T \in \Lambda}\det(CZ+D)^{-k}$$ converges uniformly on compact subsets of $M$ to a non-zero meromorphic function for each positive even integer $k>2$. Indeed $\eta_{T}(Z) = \det(CZ+D)^{n+2}$ already constitutes a complete system of factors of automorphy \cite{BG55}, and thus the functions $$z_k = \frac{\sum_{T \in \Lambda}\frac{1}{\det(CZ+D)^{(n+2)k}}}{(\sum_{T \in \Lambda}\frac{1}{\det(CZ+D)^{n+2}})^{k}}$$
generate $\invar$. Then for each $X = \sum a_k\partial_{k}$ with $a_k \in \invar$, $\partial_k = \frac{\partial}{\partial z_k}$, and $H \in \invar\Xinv$, we obtain Siegel modular equivariants $$Z_H = Z+\dot{Z}[\dot{Z}H(Z)-\half\ddot{Z}]^{-1}\dot{Z} \,\in\, \mero_\Gamma(M,\glnC).$$
\end{example}

\bibliographystyle{plain}
\bibliography{ref}
\end{document}